\documentclass[12pt,dvips]{amsart}
\usepackage{amsfonts, amssymb, latexsym, epsfig,epic,color}
\renewcommand{\familydefault}{ppl}

\newcommand\lie[1]{{\mathfrak #1}}

\newcommand\Proj{{\rm Proj\,}}
\newcommand\iso{{\cong}}
\newcommand\tensor{{\otimes}}

\newtheorem{Theorem}{Theorem} 
\newtheorem{Proposition}{Proposition} 
\newtheorem{Lemma}{Lemma}

\newtheorem*{Corollary}{Corollary}
\newtheorem*{Conjecture}{Conjecture}
 
\newtheorem*{Theorem*}{Theorem}

\newcommand\union{\bigcup}

\newcommand\CP{{\mathbb C \mathbb P}}

\newcommand\reals{{\mathbb R}}
\newcommand\complexes{{\mathbb C}}
\newcommand\integers{{\mathbb Z}}

\newcommand\naturals{{\mathbb N}}

\newcommand\Grrn{{\rm Gr}_r(\complexes^n)}
\newcommand\Gr{{\rm Gr}}
\newcommand\Hom{{\rm Hom}}

\theoremstyle{plain}

\newcommand\dfn{\bf} % maybe should be \em

\newcommand\HONEY{{\tt HONEY}}
\newcommand\BDRY{{\tt BDRY}}
\newcommand\GLn{{{GL_n(\complexes)}}}

\newcommand\LRn{\BDRY(n)}
\newcommand\ovl{\oplus}

\begin{document}
\pagestyle{plain}

\title{The honeycomb model of $\GLn$ tensor products II:\\
Puzzles determine facets of the Littlewood-Richardson cone}
\author{Allen Knutson}
\email{allenk@math.berkeley.edu}
\thanks{AK was supported by an NSF Postdoctoral Fellowship, 
  an NSF grant, and the Clay Mathematics Institute.}
\address{Mathematics Department\\ UC Berkeley\\ Berkeley, California}
\author{Terence Tao}
\thanks{TT was supported by the Clay Mathematics Institute,
  and grants from the Sloan and Packard Foundations.} 
\email{tao@math.ucla.edu}
\address{Mathematics Department\\ UCLA\\ Los Angeles, California}
\author{Christopher Woodward}
\thanks{CW was partially supported by an NSF Postdoctoral Fellowship, and
NSF Grant 9971357.}
\email{ctw@math.rutgers.edu}
%\address{Mathematics - Hill Center, 110 Frelinghuysen Rd, Piscataway NJ 08854}
\address{Mathematics Department\\ Rutgers University\\ New Brunswick, New Jersey}
\date{\today}

\maketitle

\begin{abstract}
  The set of possible spectra $(\lambda,\mu,\nu)$ of zero-sum triples
  of Hermitian matrices forms a polyhedral cone \cite{H}, whose facets
  have been already studied in \cite{Kly,HR,T,Be} in terms of Schubert
  calculus on Grassmannians. We give a complete determination of these
  facets; there is one for each triple of Grassmannian Schubert cycles
  intersecting in a unique point. In particular, the list of inequalities 
  determined in \cite{Be} to be sufficient is in fact minimal. 

% (the remaining problem was to establish the converse of the
%  result in \cite{Be}).
  
  We introduce {\em puzzles}, which are new combinatorial gadgets to compute
  Grassmannian Schubert calculus, and seem to have much % combinatorial 
  interest in their own right. As the proofs herein indicate,
  the Hermitian sum problem is very naturally studied using puzzles directly,
  and their connection to Schubert calculus is quite incidental to our
  approach.  In particular, we get new, puzzle-theoretic, proofs of
  the results in \cite{H,Kly,HR,T,Be}.

  Along the way we give a characterization of ``rigid'' puzzles, 
  which we use to prove a conjecture of W. Fulton: 
  ``if for a triple of dominant weights 
  $\lambda,\mu,\nu$ of $\GLn$ the irreducible representation 
  $V_\nu$ appears exactly once in $V_\lambda \tensor V_\mu$,
  then for all $N\in \naturals$, $V_{N\lambda}$ appears exactly once
  in $V_{N\lambda}\tensor V_{N\mu}$.''
\end{abstract}

\tableofcontents

\section{Introduction, and summary of results}

We continue from \cite{Hon1} the study of the cone
$\LRn$, which is the set of triples of weakly decreasing $n$-tuples
$(\lambda,\mu,\nu) \in (\reals^n)^3$ 
satisfying three conditions proved there to be equivalent:
\begin{enumerate}
\item regarding $\lambda,\mu,\nu$ as spectra of $n\times n$ Hermitian matrices,
  there exist three Hermitian matrices with those spectra whose
  sum is the zero matrix;
\item (if $\lambda,\mu,\nu$ are integral) regarding $\lambda,\mu,\nu$
  as dominant weights of $\GLn$, the tensor product $V_\lambda \tensor
  V_\mu \tensor V_\nu$ of the corresponding irreducible
  representations has an invariant vector;
\item regarding $\lambda,\mu,\nu$ as possible boundary data on a honeycomb,
  there exist ways to complete it to a honeycomb.
\end{enumerate}
In the present paper we
determine the minimal set of inequalities defining this cone,
along the way giving new proofs of the results in \cite{HR,T,Kly,Be}
which gave a sufficient list of inequalities in terms of
Schubert calculus on Grassmannians.
We show that this list is in fact minimal
(establishing the converse of the result in \cite{Be}).
As in \cite{Hon1}, our approach to this cone is in the honeycomb
formulation.  We also replace the use of Schubert calculus by 
{\em puzzles}, defined below.

\subsection{Prior work.}
Most prior work was stated in terms of the sum-of-Hermitian matrices problem.
That was first proved to give a polyhedral cone in \cite{H}\footnote{%
In fact Horn only proves that the cone is locally polyhedral;
convexity follows from nonabelian convexity theorems
in symplectic geometry.}.
Many necessary inequalities were found (see \cite{F1} for a survey), 
culminating in the list of Totaro \cite{T}, 
Helmke-Rosenthal \cite{HR}, and Klyachko \cite{Kly} -- hereafter we
call this the H-R/T/K result.
% ; this list is defined in terms of Schubert calculus on Grassmannians.
Klyachko proved also that this list is sufficient.
A recursively defined list of inequalities had been already 
conjectured in \cite{H};
this conjecture is true, and in fact gives the
same list as Klyachko's -- see \cite{Hon1}.

One of us (CW) observed that this list is redundant --
some of the inequalities given do not determine facets but only
lower-dimensional faces of $\LRn$ -- and proposed a criterion
for shortening the list (again in terms of Schubert calculus).
That this shorter list is already sufficient was proved by Belkale \cite{Be}.
Our primary impetus for the present work was to prove the converse:
each of these inequalities is essential, i.e. determines a facet
of $\LRn$.

\subsection{Puzzles.}
\begin{figure}[htbp]
  \begin{center}
    \leavevmode
    \epsfig{file=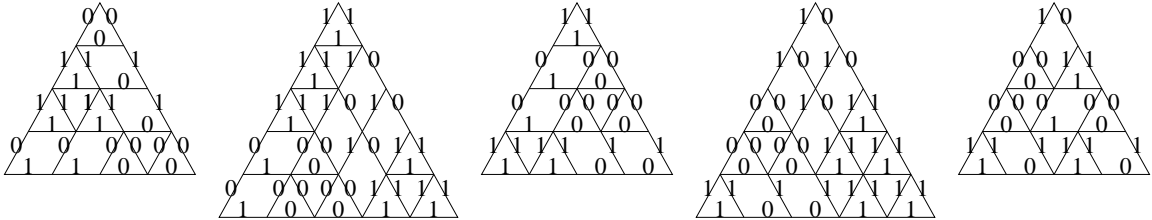,width=6in}
    \caption{Some examples of puzzles.}
    \label{fig:puzex}
  \end{center}
\end{figure}

A puzzle will be a certain kind of diagram in the triangular lattice
in the plane.  There are three {\dfn puzzle pieces}:
\begin{enumerate}
\item unit equilateral triangles with all edges labeled $0$
\item unit equilateral triangles with all edges labeled $1$
\item unit rhombi (two equilateral triangles joined together)
  with the outer edges labeled $1$ if clockwise of an obtuse angle,
  $0$ if clockwise of an acute angle.
\end{enumerate}
A {\dfn puzzle} of size $n\in \naturals$ is a decomposition of a 
lattice triangle of side-length $n$ into lattice polygons,
all edges labeled $0$ or $1$, such that each region is a puzzle piece.
Some examples are in figure \ref{fig:puzex}.

The main result about these puzzles (theorem \ref{thm:puzcount}, stated below,
proved in section \ref{sec:puzresults})
is that % suitably interpreted, 
they compute Schubert calculus on
Grassmannians. While there are many other rules for such computations,
e.g.  the Littlewood-Richardson rule, this one has the greatest
number of manifest symmetries.  (A lengthy discussion of this 
will appear in \cite{Puz2}.)

Readers only interested in a solution to the Hermitian sum problem can
skip the statement of this theorem and, in fact, 
quit after section \ref{sec:gentle}.
A central principle in the current paper is that in determining the 
facets of $\LRn$, the connection to Schubert calculus is quite irrelevant,
and it is more natural combinatorially to work with the puzzles directly,
which we do until section \ref{sec:puzresults}. This
nicely complements the principle of \cite{Hon1}, in which we worked
not with triples of Hermitian matrices but used honeycombs as their
combinatorial replacement.  We will not in general take space to repeat
the honeycomb-related definitions from \cite{Hon1}.

We fix first our conventions to describe ``Schubert calculus,''
which in modern terms is the ring structure on the cohomology
of Grassmannians.
To an $n$-tuple $\sigma$ like $00101\ldots 110$ of $r$ ones and $n-r$ zeroes, 
let $\complexes^\sigma$ denote the corresponding coordinate $r$-plane
in $\complexes^n$, and $X_\sigma$ the {\dfn Schubert cycle} 
defined as 
$$ \big\{ V_r \in \Grrn \quad: \quad
        \dim (V_r \cap F_i) \geq \dim (\complexes^\sigma \cap F_i),\quad
        \forall i\in [1,n] \big\} $$
where $\{F_i\}$ is the standard flag in $\complexes^{n}$.  
Alternately, $X_\sigma$ is the closure of the set of $r$-subspaces 
$V_r \leq \complexes^{n}$ such that 
$\sigma_i = \dim ((V_r \cap F_i) / (V_r \cap F_{i-1})), i\in [1,n]$.
The {\dfn Schubert class} $S_\sigma \in H^*(\Grrn)$ 
is the Poincar\'e dual of this cycle.
%In particular the degree of $S_\sigma$ is twice the number of pairs $i<j$
%with $1 = \sigma_i > \sigma_j = 0$.
These are well-known to give a basis for the cohomology ring.

\begin{Theorem}\label{thm:puzcount}
  Let $\pi,\rho,\sigma$ be three $n$-tuples of $r$ ones and $n-r$ zeroes,
  indexing Schubert classes $S_\pi,S_\rho,S_\sigma$ in $H^*(\Grrn)$.
  Then the following (equivalent) statements hold:
  \begin{enumerate}
  \item The intersection number $\int_{\Grrn} S_\pi S_\rho S_\sigma$ is
    equal to the number of puzzles whose
    NW boundary edges are labeled $\pi$, NE are labeled $\rho$, and
    S are labeled $\sigma$, all read clockwise.
%  \item The structure constant $c^\sigma_{\pi\rho}$ is
%    equal to the number of puzzles whose
%    NW boundary edges are labeled $\pi$, NE are labeled $\rho$, and
%    S are labeled $\sigma$, all read left-to-right.
  \item $$ S_\pi\, S_\rho 
  = \sum_{\hbox{puzzles $P$ } }
  S_{\hbox{southern side of $P$, read left to right}} $$ 
        where the sum is taken over puzzles with NW side labeled $\pi$,
        NE side labeled $\rho$, both from left to right.
  \end{enumerate}
\end{Theorem}

\begin{figure}[htbp]
  \begin{center}
    \leavevmode
    \epsfig{file=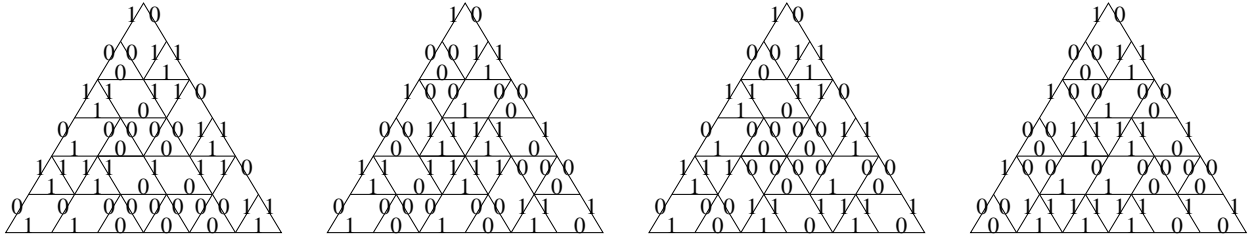,width=6.5in}
    \caption{The four puzzles $P$ with NW and NE boundaries each labeled
      $010101$, computing $S_{010101}^2 = S_{110001} + 2 S_{101010} +
      S_{011100}$ in $H^*(\Gr_3(\complexes^6))$.}
%{\bf check the Schubert polynomial product for extra, unstable terms}
    \label{fig:prodex}
  \end{center}
\end{figure}

This first is the advertised $Z_3$-invariant formulation. The second
formulation is very suitable for computations; an example is in figure
\ref{fig:prodex}.

Puzzles have another symmetry, which we call {\dfn puzzle duality}: 
the {\dfn dualization} of a puzzle is defined to
be the left-right mirror reflection, with all $0$s exchanged for $1$s
and vice versa.
This realizes combinatorially another symmetry of Schubert calculus,
coming from the isomorphism of the $r$-Grassmannian in an 
$n$-dimensional space $V$ with the $(n-r)$-Grassmannian in $V^*$.
We will use puzzle duality to reduce the number of cases considered
in some arguments. % to come.

\subsection{Organization of this paper.}
% This paper is organized as follows.
In sections \ref{sec:hr}-\ref{sec:gentle} we classify the facets of
$\LRn$ in terms of puzzles. In sections \ref{sec:puzresults}-\ref{sec:schub}
we prove and make use of the connection of puzzles to Schubert calculus.
To emphasize again: the reader who is only looking for the minimal
list of inequalities determining $\LRn$ may completely ignore this
connection, and take puzzles as the more relevant concept than
Schubert calculus! 

Here is a slightly more detailed breakdown of the paper.
In section \ref{sec:hr} we prove the puzzle-theoretic analogue
of the H-R/T/K result: each puzzle gives an
inequality on $\LRn$.

In section \ref{sec:tco} we essentially repeat Horn's analysis of
the facets of $\LRn$, but in the honeycomb framework; the analogues
of his direct sums turn out to be {\em clockwise overlays}.
Using these we prove the puzzle-theoretic analogue of Klyachko's
sufficiency result (a converse of H-R/T/K): 
every facet comes from a puzzle. Easy properties of puzzles
(from section \ref{sec:puzresults}) then imply Horn's results
(but not his conjecture).

In section \ref{sec:gentle} we study ``gentle loops'' in puzzles, 
and show that the minimal list of inequalities is given by 
puzzles with no gentle loops. Then comes the only particularly
technical part of the paper: showing that puzzles without gentle loops
are exactly the {\em rigid} ones, meaning those determined by their 
boundary conditions. That the rigid-puzzle inequalities are a 
sufficient list is the puzzle analogue of Belkale's result \cite{Be};
conversely,
that every rigid-puzzle inequality determines a facet, is the central
new result of this paper.

In section \ref{sec:puzresults} we describe the connection of puzzles 
to Schubert calculus, and in section \ref{sec:schub}
give puzzle-free statements of our theorems. This section also
serves as a summary of the old and new results in this paper.

Since Schubert calculus is itself related
to the tensor product problem (in a lower dimension), this gives a
combinatorial way to understand the still-mysterious Horn recursion.
We also give an application of the no-gentle-loop characterization
of rigid puzzles to prove an unpublished conjecture of W. Fulton.

In the last section we state the corresponding results for sums of $m$
Hermitian matrices.  The proofs extend almost without change to the $m\geq 3$
case. 
In an appendix we give a quick proof of the equivalence between
the three definitions of $\LRn$, replacing Klyachko's argument by
the Kirwan/Kempf-Ness theorem, which allows for rather stronger results.

Since completing this work, we received the preprint \cite{DW1}, 
which studies representations of general quivers; our results can be
seen as concerning the very special case of the ``triple flag quiver''.
Assuming Fulton's conjecture as input (see conjecture 30 of \cite{DW1}), 
their results provide a (completely different) proof of the converse 
of Belkale's result.
We have been unable to find any generalization of our
honeycomb and puzzle machinery to general quivers.

We are most grateful to Anda Degeratu for suggesting the name
``puzzle'', and the referee for many useful comments.  The
once-itinerant first author would like to thank Rutgers, UCLA, MSRI,
and especially Dave Ben-Zvi for their gracious hospitality while part
of this work was being done.

\section{Puzzles give inequalities on $\LRn$}\label{sec:hr}

In this section we determine a list of inequalities satisfied by $\LRn$, 
which will eventually be seen to be the puzzle-theoretic version of
the H-R/T/K result.

Recall from \cite{Hon1} that $\LRn$ is defined as the image of the 
``constant coordinates of boundary edges'' 
map $\partial: \HONEY_n \to (\reals^n)^3$. In proposition 1 of \cite{Hon1}
we showed that the nondegenerate honeycombs (those whose edges are all 
multiplicity 1 and vertices all trivalent) are dense in $\HONEY_n$. 
So in determining inequalities on $\LRn$ one can safely restrict
to boundaries of nondegenerate honeycombs. This reduction is
not logically necessary for the rest of the section but may
make it easier to visualize.

% \subsection{Puzzles and inequalities on $\LRn$.}

\begin{figure}[htbp]
  \begin{center}
    \leavevmode
    \epsfig{file=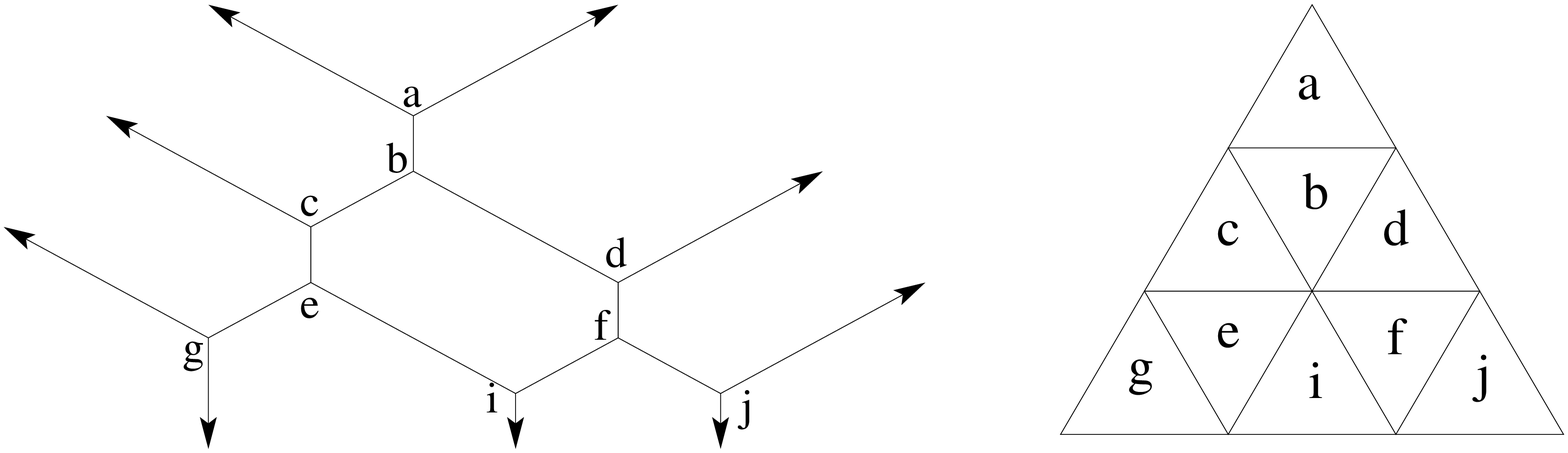,height=1.5in}
    \caption{Vertices of nondegenerate honeycombs 
      correspond to triangles in puzzles. The $n=3$ case is pictured.}
    \label{fig:puzcor}
  \end{center}
\end{figure}

Let $h$ be a nondegenerate $n$-honeycomb, and let $\pi$ be a lattice
equilateral triangle of side-length $n$.
There is an obvious correspondence between $h$'s vertices and 
the unit triangles in $\pi$, as in figure \ref{fig:puzcor}.
More importantly for us, one can also correspond the bounded edges of $h$
(connecting two vertices) and the unit rhombi in $\pi$ (the union of
two triangles). Finally,
the semiinfinite edges in $h$ correspond to
the boundary edges of $\pi$ (to which they are perpendicular).

Given a puzzle $P$ of side-length $n$, define a linear functional 
$f_P: \HONEY_n \to \reals$ by
$$ f_P(h) = \sum_{\hbox{$R$ rhombus in $P$}} 
        \hbox{ the length of the corresponding edge in $h$. } $$ 
(This does not really use $h$ nondegenerate -- in the degenerate case,
some of these terms are $0$.) Note that this functional is automatically
nonnegative, as it is a sum of nonnegative terms.
We define ``length'' relative to the triangular lattice, i.e.
$1/\sqrt{2}$ of the usual Euclidean length.
% of the line between $(0,0,0)$ and $(1,0,-1)$ is taken to be $1$.

As defined above, the quantity $f_P(h)$ seems to depend on the
internal structure of $P$ and $h$.  However, there is a ``Green's theorem''
which allows us to write $f_P(h)$ purely in terms of the boundary
labels on $P$ and $h$:

\begin{Theorem}\label{thm:hrpuz}
  Let $P$ be an $n$-puzzle, and $f_P$ the corresponding functional 
  on $\HONEY_n$. Then 
$$ f_P(h) = \sum_{\hbox{$P$'s boundary edges labeled $0$}}
   \hbox{the constant coordinate on the corresponding boundary edge of $h$} $$ 
and in particular descends to give a nonnegative functional on $\LRn$.
Put another way, the inequality $f_P\geq 0$ is satisfied by $\LRn$.
\end{Theorem}

\begin{proof}
  We compute what at first seems to be a different functional,
  in two different ways.
  Call an edge on a puzzle piece {\dfn right-side-up} if its outward normal
  is parallel to an outward normal of the entire puzzle, 
  {\dfn upside-down} if the 
  outward normal is antiparallel. So on a right-side-up triangle, all
  three edges are right-side-up, and vice versa for an upside-down triangle.
  Whereas on a rhombus, two of the edges are right-side-up, two upside-down.

  Define the functional $g_P: \HONEY_n \to \reals$ by
$$ g_P(h) = \sum_{\hbox{$p$ a piece}} \,\,\sum_{\hbox{$e$ a $0$-edge of $p$}}
(\pm 1)
\hbox{constant coordinate on the corresponding edge of $h$} $$ 
where the sign is $+1$ if $e$ right-side-up, $-1$ if $e$ upside-down.

We claim first that $g_P = f_P$. Consider the contribution a piece $p$
makes to the sum in $g_P$: a $1,1,1$-triangle contributes nothing,
a $0,0,0$-triangle contributes the three coordinates of a vertex
(which sum to zero), and we leave the reader to confirm that a
rhombus contributes the length of the corresponding edge in $h$.

To show that $g_P$ also matches the conclusion of the theorem,
rewrite by switching the order of summation:
$$ g_P(h) = \sum_{\hbox{$e$ a $0$-edge of $P$}} \,
\sum_{\hbox{$p$ containing $e$}} (\pm 1)
\hbox{constant coordinate on the corresponding edge of $h$} $$ 
$$  = \sum_{\hbox{$e$ a $0$-edge of $P$}} 
\hbox{constant coordinate on the corresponding edge of $h$} 
\sum_{\hbox{$p$ containing $e$}} (\pm 1) $$ 
For every edge $e$ internal to the puzzle, this latter sum is $+1-1$
which cancels, whereas for every exterior $0$-edge it is $1$.
The claim follows.

And as stated before, this functional is a sum of honeycomb edge-lengths, 
so automatically nonnegative on $\HONEY_n$.
\end{proof}

\begin{figure}[Hht]
  \begin{center}
    \leavevmode
    \input{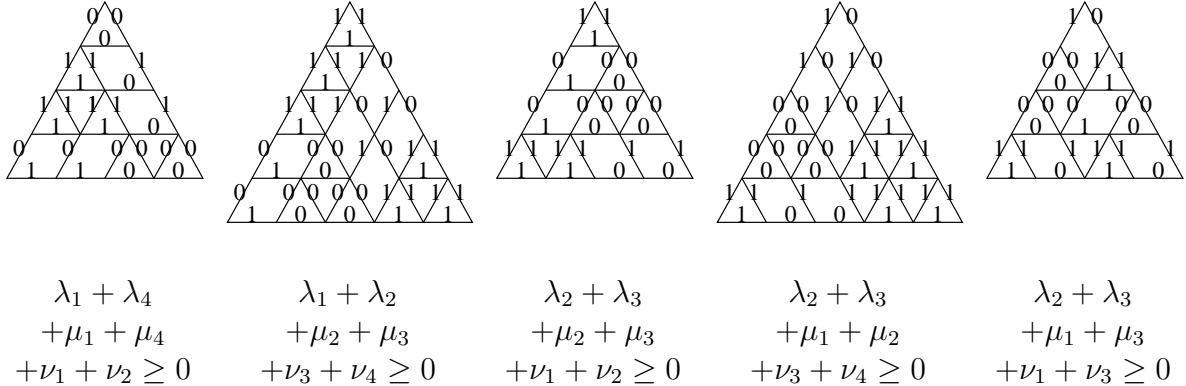}
    \caption{Some puzzles and the corresponding inequalities on $\LRn$.
      The subscripts on $\lambda$ correspond to the locations of the $0$s
      on the NW side, $\mu$ on the NE, and $\nu$ on the South.}
    \label{fig:ineqex}
  \end{center}
\end{figure}

(As we will review in subsection \ref{ss:horn},
these inequalities $f_P \geq 0$ are automatically of the sort that
Horn predicted in \cite{H} -- a sum of $r$ distinct
elements from each of $\lambda$, $\mu$, and $\nu$.)
In figure \ref{fig:ineqex} we repeat the puzzles from figure \ref{fig:puzex}
and give the corresponding inequalities.

Since the sum of {\em all} the boundary coordinates is zero, this inequality
can be restated as coming from a nonpositive functional,
$$ \sum_{\hbox{ $e$ a boundary edge of P}} \hbox{(label on $e$)} \cdot
  \hbox{(constant coordinate on $h$'s edge corresponding to $e$)} \leq 0, $$ 
as it does in some of the literature (e.g. \cite{F1}).

This theorem \ref{thm:hrpuz} will in section \ref{sec:schub} be seen
to be the puzzle analogue of the necessary conditions of
H-R/T/K.
In the next section we show that every facet (except for some easy,
uninteresting ones) does indeed come from a puzzle inequality $f_P \geq 0$, 
which will be the puzzle-theoretic analogue of
Klyachko's sufficiency theorem (a converse of H-R/T/K).

\section{Facets come from puzzles, via clockwise overlays}\label{sec:tco}

In this section we study the honeycombs that lie over facets of $\LRn$.
We begin by recalling Horn's results \cite{H} on triples of Hermitian
matrices which sum to zero, to help make intuitive the corresponding
results we will find on honeycombs.

\subsection{Horn's results.}\label{ss:horn}
Horn considered the function ``take eigenvalues in decreasing order'' 
from zero-sum Hermitian triples to $(\reals^n)^3 = \{(\lambda,\mu,\nu)\}$. 
By definition, the image satisfies the {\dfn chamber inequalities}, i.e. that
$\lambda_i \geq \lambda_{i+1}\,\forall i$ (similarly $\mu,\nu$). 
The remaining facets we call {\dfn regular facets}.

Away from the chamber walls, the ``take eigenvalues'' map is differentiable,
and one can use calculus to find its extrema:
Horn did this, and found that the critical points occur exactly when
the zero-sum Hermitian triple is a direct sum of two smaller ones.
(This is nowadays a standard calculation in Hamiltonian geometry --
see \cite{K} for an exposition of this viewpoint.)
Since that implies that the traces of each subtriple sum to zero,
one sees that the equation of the facet so determined says that
the sum of a certain $r$ eigenvalues from $\lambda$, 
another $r$ from $\mu$, and another $r$ from $\nu$ add to zero.

Also, the Hessian is definite at an extremal point, which gives another
condition on these three $r$-element subsets $I,J,K$ of $\{1,\ldots,n\}$. 
Define an {\dfn inversion} of such a subset $S$ as  
a pair $(a<b) \subseteq \{1,\ldots,n\}$ such that $a\in S, b\notin S$.
Then Horn shows that definiteness of the Hessian implies that the total
number of inversions, over the three subsets $I,J,K,$ is $r(n-r)$.
(Both of these conditions are automatic for puzzles, as shown later in
proposition \ref{prop:puzstructure}; in particular this will give
combinatorial proofs of Horn's results.)

We now undertake the same extremal analysis on honeycombs, rather than
zero-sum Hermitian triples. We will need the following lemma, whose proof 
is immediate, to recognize inequalities from individual boundary points.
Recall that a {\dfn facet} of a polyhedron is a codimension-1 face.

\begin{Lemma}\label{lem:bdry}
  Let $P$ be a polyhedron (convex, but not necessarily compact), 
  $p$ a point on a facet $\Phi$ of $P$,
  and $f$ a nonzero affine-linear function vanishing at $p$. 
  If $\Phi$ contains a neighborhood of $p$ in $f^{-1}(0)$, 
  then $p$ is an interior point of $\Phi$, 
  the equation of $\Phi$ is $f = 0$, 
  and the inequality determining $\Phi$ is either $f\geq 0$ or $f\leq 0$.
\end{Lemma}

To apply this to $\LRn$, we will need to know its dimension;
via the Hermitian picture, this is well known to be $3n-1$
(it is cut down from $3n$ by the fact that the sum of the
traces must be zero). We give a honeycomb-theoretic proof in 
proposition \ref{prop:LRndim},
mainly in order to introduce the construction by which we will
vary the boundary of a honeycomb.

Define the {\dfn natural sign} of an oriented edge in a honeycomb to
be $+1$ if the edge points Northwest, Northeast, or South and $-1$ if
it points North, Southwest, or Southeast.  (By these six compass
directions we of course really mean directions that are at $60^\circ$
angles from one another, not $45^\circ$ and $90^\circ$.) Observe that
a path in a nondegenerate honeycomb must alternate natural sign
(orienting the edges to follow the path). In particular, a path coming
in from infinity on one boundary edge (natural sign $-1$) and going
out on another (natural sign $+1$) must be of odd length.

\begin{Proposition}\label{prop:LRndim}
  The cone $\LRn$ is $(3n-1)$-dimensional.
\end{Proposition}

\begin{proof}
  This is certainly an upper bound: by lemma 1 of \cite{Hon1}, 
  the sum of all the constant coordinates of boundary edges is zero. 
%  (In the Hermitian 
%  picture, this is a trace condition; in the representation theory picture,
%  this says that the center acts with weight zero on invariant vectors.)
  
  Let $h$ be a nondegenerate honeycomb, $\epsilon$ a (possibly
  negative) real number such that $2|\epsilon|$ is smaller than the
  length of any of $h$'s edges, and $e,f$ two boundary edges.
  Then there exists a path $\gamma$ in the honeycomb 
  tinkertoy $\tau_n$ connecting $e$ and $f$ (which we can ask be
  non-self-intersecting).
  We can add $\epsilon$ times the natural sign to the constant coordinates
  of $h$'s edges along $\gamma$ and get a new honeycomb.
  
  This changes $e$'s coordinate by $-\epsilon$, and $f$'s by
  $+\epsilon$. By repeating this with other pairs, we can achieve
  arbitrary small perturbations of the boundary coordinates, subject
  to the sum staying zero. So $\LRn$ contains a $(3n-1)$-dimensional
  neighborhood of $\partial h$.
\end{proof}

Call the construction in proposition \ref{prop:LRndim} the
{\dfn trading construction}. We will need it not only for the
honeycomb tinkertoy $\tau_n$, but (connected) tinkertoys constructed from 
$\tau_n$ by eliding simple degeneracies, as we did
in the corollary to theorem 1 of \cite{Hon1}. 
In particular if $h$ is a simply degenerate honeycomb, and $\tau_n$
stays connected after eliding $h$'s simple degeneracies, then
$\partial h$ is in the interior of $\LRn$.
% The oriented paths in the trading construction have
% an additional interesting property: if one labels an oriented edge
% with an {\dfn outerwardness sign} of $+1$ if it is parallel to an
% outward-pointing unbounded edge, $-1$ if antiparallel, then the signs
% assigned to edges during the trading construction in proposition
% \ref{prop:LRndim} either all match or all anti-match their 
% outwardness signs. 

\subsection{Extremal honeycombs are clockwise overlays.}
Recall the overlay operation from \cite{Hon1}; it makes an
$n$-honeycomb $A\ovl B$ from an $(n-r)$-honeycomb $A$ and an
$r$-honeycomb $B$.
If $p$ is a point common to $A$ and $B$, call it a {\dfn transverse
point of intersection} if it is a vertex of neither, and isolated
in the intersection. Call $A\ovl B$ a {\dfn transverse overlay} 
if all intersection points are transverse. In this case every small
perturbation of $A$ and $B$ is again a transverse overlay; by proposition
\ref{prop:LRndim} this gives a $(3r-1) + (3(n-r)-1) = 3n-2$ dimensional
family of boundaries.

\begin{figure}[htbp]
  \begin{center}
    \leavevmode
    \epsfig{file=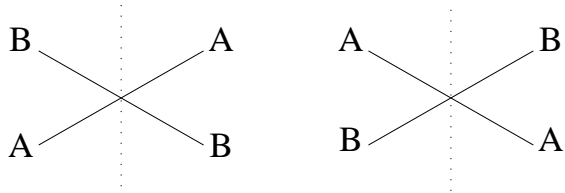,height=1in}
    \caption{In the left figure, $A$ turns clockwise to $B$, whereas in
      the right the reverse is true. Any transverse point of intersection
      of two overlaid honeycombs must look like exactly one of these,
      up to rotation.}
    \label{fig:clockwise}
  \end{center}
\end{figure}

If $p$ is a transverse intersection point of $A\ovl B$, then up to rotation
a neighborhood of $p$ looks like exactly one of the two pictures
in figure \ref{fig:clockwise}; say that {\dfn $A$ turns clockwise to $B$ 
at $p$} or {\dfn $B$ turns clockwise to $A$ at $p$} depending on which.

\begin{figure}[htbp]
  \begin{center}
    \leavevmode
    \epsfig{file=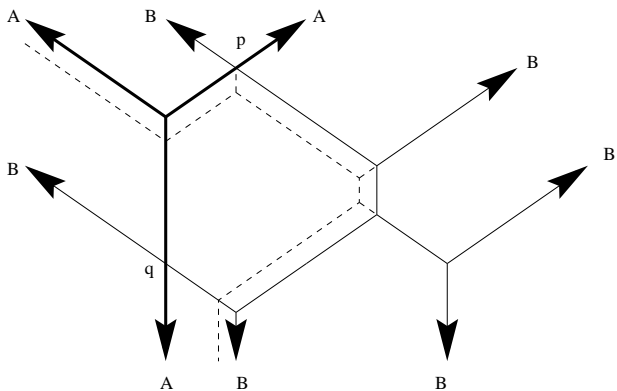,height=2in}
    \caption{The solid honeycomb $A$ is clockwise to the thin one $B$
      at the point $p$, vice versa at $q$.
      The dashed line indicates the result of trading an $A$ edge up
      and a $B$ edge down, using a path through the crossing $p$.}
    \label{fig:breatheclock}
  \end{center}
\end{figure}

Recall from \cite{Hon1} that a simple degeneracy of a honeycomb is a
vertex where two multiplicity-one edges cross in an X.
We can deform a simple degeneracy to a pair of nondegenerate vertices
connected by an edge. If we do this in an overlay as in figure 
\ref{fig:breatheclock},
who is clockwise to whom determines the 
resulting behavior of the boundary, as explained in the following lemma.

\begin{Lemma}\label{lem:clockwisetrade}
  Let $h = A \ovl B$ be a transverse overlay of two nondegenerate
  honeycombs, and $p$ a point of intersection, such that $A$ turns
  clockwise to $B$ at $p$.  Let $\gamma_A$ be a path in $A$ that comes from
  infinity to $p$, and $\gamma_B$ a path in $B$ that goes from $p$
  to infinity. Then we can extend the trading construction to 
  $\gamma_A \union \gamma_B$ and increase the constant coordinate
  on the first edge of $\gamma_A$ while decreasing the constant
  coordinate on the last edge of $\gamma_B$ the same amount,
  leaving other boundary edges unchanged. 
\end{Lemma}

\begin{proof}
  By rotating if need be, we can assume $p$ looks like the left figure
  in figure \ref{fig:clockwise}, a simple degeneracy. Assume the path
  $\gamma_A$ comes from the Southwest, $\gamma_B$ going to the Southeast
  (the other three cases are similar).
  
  We can pull the edges in $\gamma_A,\gamma_B$ at $p$ down and create
  a vertical edge in the middle (as in the figure \ref{fig:breatheclock}
  example). To extend this change to the rest of the 
  honeycomb, we apply the trading construction to $\gamma_A\union\gamma_B$,
  but we must move edges a {\em negative} $\epsilon$ times their
  natural signs (or else the vertical edge created will have negative
  length). In particular the first edge of $\gamma_A$, whose natural
  sign is negative, has its constant coordinate increased.
\end{proof}

The definition of ``largest lift with respect to a superharmonic functional''
was one of the more technical ones from \cite{Hon1}; the details of it
are not too important in the following lemma, except for the
application of \cite{Hon1}'s theorem 2 (as explained within).
  
\begin{Lemma}\label{lem:llco}
  Let $b$ be a generic point on a regular facet $\Phi$ of $\LRn$,
  and $h$ a largest lift of $b$ (with respect to some choice of
  superharmonic functional on $\HONEY_n$). 
  Then $h$ is a transverse overlay $A\ovl B$ of two smaller honeycombs,
  where at every point $p$ of intersection $A$ turns clockwise to $B$.

%  If $b$ is in the interior of the facet $\Phi$,
  In addition,
  the inequality determining $\Phi$ says that the sum of the
  constant coordinates of $h$'s boundary edges contained in $A$ 
  is nonnegative.
\end{Lemma}

The genericity condition on $b$ is slightly technical: 
we ask that at most one proper subset of the boundary edges 
(up to complementation) has total sum of the constant coordinates being zero. 
This avoids a finite number of $(3n-3)$-dimensional subspaces of the
$(3n-2)$-dimensional facet $\Phi$, and as such is an open dense condition.
Also we ask that $b$ be regular.

\begin{proof}
  By theorem 2 of \cite{Hon1}, $h$ is simply degenerate and acyclic
  (this uses $b$ regular and $h$ a largest lift).
  Let $\tau$ be the post-elision tinkertoy, of which $h$
  can be regarded as a nondegenerate configuration. We first claim that
  $\tau$ is disconnected. For otherwise, we could use the trading
  construction to vary the boundary of $h$ in arbitrary directions
  (subject to the sum of the coordinates being zero), and therefore
  $b$ would not be on a facet of $\LRn$.
  So we can write $\tau = \rho \cup \sigma$, and $h = A\ovl B$,
  where $A,B$ are honeycombs with tinkertoys $\rho,\sigma$ respectively.
  
  The boundary $b = \partial h$ therefore satisfies the equation 
  ``the sum of the boundary coordinates of $h$ belonging to $\rho$ is zero.''
  (Likewise $\sigma$.) This remains true if we deform $A$ and $B$
  as individual honeycombs. By the genericity condition on $b$,
  $A$ and $B$ must each be connected, so varying them 
  gives us a $(3n-2)$-dimensional family of variations of $b$ with
  $b$ in the interior.
% (If $A$ is of size $k$, $B$ of size $n-k$,
%  varying them independently gives families of dimensions $3k-1$
%  and $3(n-k)-1$, for a total of $3n-2$.)  
  By lemma \ref{lem:bdry} we have found the equation for 
  the facet containing $b$.
  
  It remains to show that $A$ turns clockwise to $B$ at every
  intersection; actually we will only show that all the intersections
  are consistent (and switch the names of $A$ and $B$ if we chose
  them wrongly). For each intersection $p$, let $\tau_p$ be the
  tinkertoy made from the honeycomb tinkertoy $\tau_n$ by 
  eliding all of $h$'s simple degeneracies {\em other than} $p$. 
  Since the fully elided $\tau$ is acyclic with two components, 
  each $\tau_p$ is acyclic with one component.

  We can now attempt to trade $A$'s boundary coordinates for $B$'s.
  Fix a semiinfinite edge of $A$ and one of $B$.
  Since $\tau_p$ is acyclic there will be only one path $\gamma$ 
  connecting them, necessarily going through $p$. 
  By lemma \ref{lem:clockwisetrade} we can apply the trading
  construction to $\gamma$, increasing the constant coordinate on 
  our semiinfinite edge of $A$
  if $A$ turns clockwise to $B$ at $p$, decreasing it in the other case.

  If there exist vertices $p,q$ such that $A$ turns clockwise to $B$ at
  $p$, but vice versa at $q$, then by trading we can move to either
  side of the hyperplane determining the facet. This contradiction
  shows that the intersections must be consistently all clockwise
  or all counterclockwise.
\end{proof}

The analogy between the Hermitian direct sum operation and the
honeycomb overlay operation is even tighter than this:
at the critical values of ``take eigenvalues'' that are not at
extrema, one also finds transverse overlays, 
and the index of the Hessian can be computed from 
the number of intersections that {\em are} clockwise. 
However, until a tighter connection is found someday in the form of,
say, a measure-preserving map from zero-sum Hermitian triples to the
polytope of honeycombs, the Hermitian and honeycomb theorems will have
to be proven independently.

This lemma \ref{lem:llco} motivates the following definition: 
say that an overlay $A\ovl B$ is a {\dfn clockwise overlay} 
(without mentioning a particular point) if the overlay is transverse, 
and at all points of intersection $A$ turns clockwise to $B$.
This is probably not the right definition:
because of the insistence on transversality, it is not closed under limits.
However, since in this paper we will only be interested in transverse overlays,
it will be more convenient to build it into the definition.

A very concrete converse to this lemma is available:

\begin{Lemma}\label{lem:witness}
  Let $h=A\ovl B$ be a clockwise overlay.
  Let $\mathcal A$ be the subset of $\tau_n$'s semiinfinite edges in
  the $A$ part of $h$.  Then the inequality
$$ \sum_{e\in \mathcal A} \hbox{the constant coordinate on $e$} \geq 0 $$ 
  defines a regular facet of $\LRn$, containing $\partial h$. Moreover,
%  this inequality is of the form $f_P \geq 0$ for some puzzle $P$
  there exists a puzzle $P$ such that this inequality is the one 
  $f_P \geq 0$ associated by theorem \ref{thm:hrpuz}.
\end{Lemma}

\begin{proof}
  Plainly $h$ satisfies this inequality with equality.  
  We can deform $A$ and $B$ to nondegenerate honeycombs $A'$ and $B'$;
  if we move the vertices of each little enough, they will not cross over
  edges of the other, and the result will again be a clockwise overlay $h'$,
  satisfying the same equality.

  Build a puzzle $P$ from $h'$ as follows: 
  \begin{itemize}
  \item to each vertex in $A'$, associate a {\dfn 0,0,0}-triangle
  \item to each vertex in $B'$, associate a {\dfn 1,1,1}-triangle
  \item to each crossing vertex in $h'$, associate a rhombus
  \end{itemize}
  with the puzzle pieces glued together if the vertices in $h$ share
  an edge. Then the fact that $A'$ turns clockwise to $B'$ means
  that the labels on the rhombi will match the labels on the triangles,
  so $P$ will be a puzzle. An example is in figure \ref{fig:tco2puz}.
  \begin{figure}[htbp]
    \begin{center}
      \leavevmode
      \epsfig{file=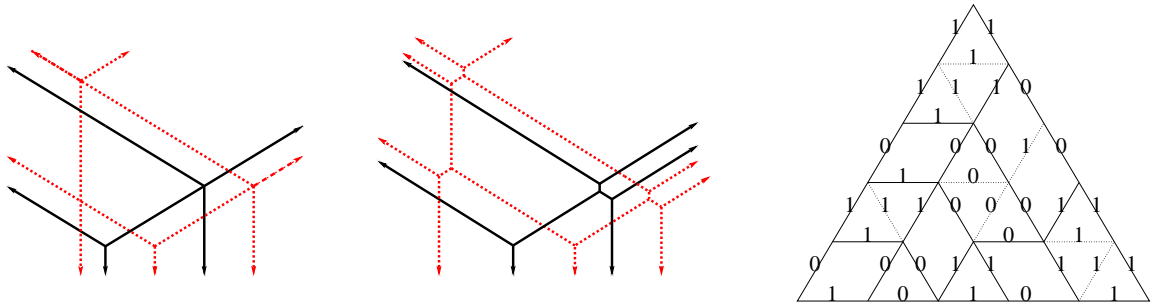,width=6in}
      \caption{A clockwise overlay, a deformation of its 
        constituents to nondegenerate honeycombs, and the puzzle built
        from that.  The corresponding inequality on $\BDRY(5)$ is
        $\lambda_1+\lambda_3+\mu_2+\mu_3+\nu_2+\nu_4\geq 0$.
        Note that some edges of the original overlay have multiplicity 2,
        leading to repetitions along the boundary of the puzzle.}
      \label{fig:tco2puz}
    \end{center}
  \end{figure}
  
  We can perturb $A'$ and $B'$ small amounts and vary the 
  boundary coordinates of in arbitrary directions. 
  This gives us a $(3n-2)$-dimensional family
  of possible boundaries containing our original point $\partial h$,
  all satisfying the stated inequality.
  By lemma \ref{lem:bdry}, $h$ is in the interior of a facet
  determined by this inequality, which is the inequality $f_P \geq 0$.
\end{proof}

Call a clockwise overlay a {\dfn witness} to the facet it exhibits
via lemma \ref{lem:witness}. This gives us a convenient way to
exhibit facets.

% (Note that from a non-clockwise overlay one could construct a ``not-quite
% puzzle,'' with some rhombi with 1s and 0s in the wrong places. 
% This is analogous to looking for critical points of the ``take
% eigenvalues'' map but not demanding them to be extremal.)

Applying lemma \ref{lem:witness} to lemma \ref{lem:llco}, we get

\begin{Theorem}\label{thm:klypuz}
  Let $\Phi$ be a regular facet of $\LRn$. 
  Then there exists a puzzle $P$ such that $\Phi$ is the facet determined
  by $f_P \geq 0$, i.e. $\Phi = \LRn \cap f_P^{-1}(0)$.
\end{Theorem}

So each regular facet gives a puzzle, and each puzzle gives an inequality,
that together with the list of chamber inequalities, 
determine $\LRn$ (and as we will see in section \ref{sec:schub}, 
this list is the same as Klyachko's, itself the same as Horn's).

But not every inequality is satisfied with equality on a facet.
Define an inequality $f \geq 0$ on a polyhedron $\Pi$ to be 
{\dfn essential} if $f^{-1}(0) \cap \Pi$ is a facet of $\Pi$, and
{\dfn inessential} if $f^{-1}(0) \cap \Pi$ is lower-dimensional.
-- equivalently, some positive multiple
of it must show up in any finite list of inequalities that determine $\Pi$.
For example, for a point in the plane to be in the first quadrant
it is necessary and sufficient that it satisfy the inequalities 
$\{x\geq 0, y\geq 0, x+y\geq 0\}$,
but the third inequality is only pressed at the origin, and
can be omitted from the list).
In the next section we will cut our list of
inequalities on $\LRn$ down to the essential inequalities.

\subsection{Independence of the chamber inequalities for $n>2$.}
We have thus far %essentially
ignored the chamber inequalities 
$ \lambda_i\geq \lambda_{i+1}$ etc. on $\LRn$, 
focusing attention on the inequalities determining regular facets.
We thank Anders Buch for pointing out to us the following subtlety
that this perspective misses.

\begin{Theorem}\label{thm:anders}
  For $n>2$, the chamber inequalities on $\LRn$ are essential.
  For $n=2$, they are implied by the regular inequalities and the equality 
  $\lambda_1 + \lambda_2 + \mu_1 + \mu_2 + \nu_1 + \nu_2 = 0$.
\end{Theorem}

\begin{proof}
  Consider a honeycomb $h$ satisfying 
  $\lambda_i = \lambda_{i+1}$ for some $i$, but otherwise
  minimally degenerate; an example is in figure
  \ref{fig:anders}. (By $Z_3$ symmetry it is enough to consider 
  the $\lambda$ case.)
  \begin{figure}[htbp]
    \begin{center}
      \leavevmode
      \epsfig{file=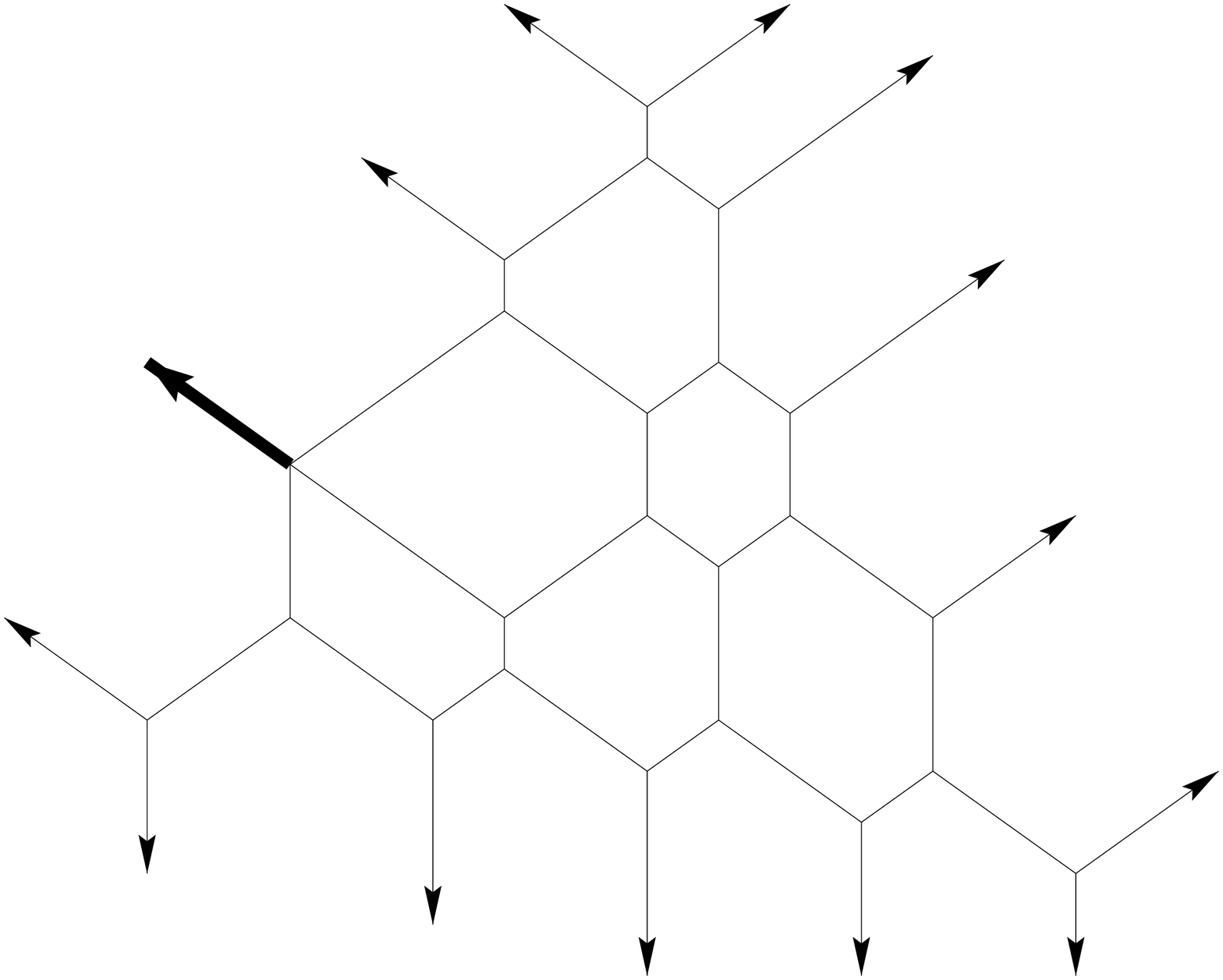,height=1.5in}
  \caption{A generic honeycomb over the chamber facet 
    determined by $\lambda_2 = \lambda_3$.}
      \label{fig:anders}
    \end{center}
  \end{figure}
  Note that there is only one nongeneric vertex, at the end of the
  multiplicity-two semiinfinite edge. It is straightforward to construct
  a similar such honeycomb for any $n$ and $i$.
  
  We mimic the proof of proposition \ref{prop:LRndim}, in using the
  trading construction to exhibit a $(3n-2)$-dimensional family in
  $\LRn$ satisfying $\lambda_i = \lambda_{i+1}$. We can move the
  doubled edge wherever we like by translating the whole honeycomb.
  To trade the constant coordinates of any other two 
  boundary edges, we use paths avoiding the bad vertex in $h$,
  which exist for $n>2$. This proves the first claim.
  
  For $n=2$ we hit a snag -- sometimes the only paths will go through
  the bad vertex. But we can check $n=2$ directly. The regular
  inequalities are
$$ \lambda_1 + \mu_1 + \nu_2 \geq 0, \quad
\mu_1 + \nu_1 + \lambda_2 \geq 0, \quad
\nu_1 + \lambda_1 + \mu_2 \geq 0. $$
  Sum the first two, and subtract the equality
  $\lambda_1 + \lambda_2 + \mu_1 + \mu_2 + \nu_1 + \nu_2 = 0$ to get
  $\mu_1 - \mu_2 \geq 0$. The other two inequalities are proved in ways 
  symmetric to this one.
\end{proof}

\section{Gentle loops vs. rigid puzzles}\label{sec:gentle}

We have at this point an overcomplete list of inequalities, 
coming from puzzles; in section \ref{sec:schub} we will see 
it is exactly that of H-R/T/K 
(which by \cite{Hon1} is exactly that of Horn's conjecture \cite{H}).
Our remaining goal is to cut down the list of inequalities to 
the essential set -- those that determine facets of $\LRn$,
rather than lower-dimensional faces.

To do this, we will shortly introduce the concept
of a {\em gentle loop} in a puzzle, and prove two things:
\begin{itemize}
\item regular facets correspond 1:1 to puzzles without gentle loops 
\item a puzzle has no gentle loops if and only if it is {\dfn rigid}, 
  i.e. is uniquely determined by its boundary conditions.
\end{itemize}
The first is remarkably straightforward, the second a bit more technical.

\begin{figure}[htbp]
  \begin{center}
    \leavevmode
    \epsfig{file=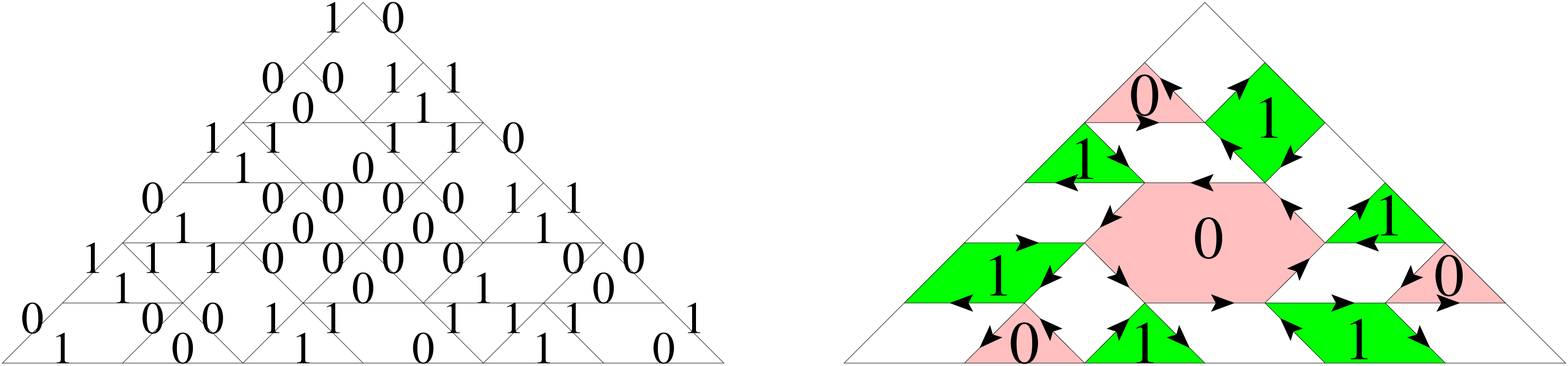,height=1.5in,width=3.5in}
    \caption{A puzzle, and its decomposition into regions, with region edges
      oriented. This is one of the two smallest examples with a gentle loop;
      it goes counterclockwise around the central hexagon.}
    \label{fig:puzregions}
  \end{center}
\end{figure}

Cut a puzzle up along the interior edges that separate two distinct types
of puzzle pieces; call the connected components of what remains the
{\dfn puzzle regions}, coming in the three types {\dfn $0$-region},
{\dfn $1$-region,} and {\dfn rhombus region}. 
Define a {\dfn region edge} in a puzzle as one separating two distinct
types of puzzle piece.  Thus every region edge either separates a
rhombus region from a $0$-region, or a $1$-region from a rhombus
region.  Orient these edges, so that $0$-regions are always on the
left, and $1$-regions are always on the right. (Viewed as edges of the
parallelograms, this orients them to point away from the acute vertices.
Stated yet another way, 
they go clockwise around the $1$-regions, counterclockwise
around the $0$-regions.) An example is in figure \ref{fig:puzregions}.

Define a {\dfn gentle path} in a puzzle as a finite list of region edges,
such that the head of each connects to the tail of the next, and
the angle {\em of turn} is either $0^\circ$ or $60^\circ$, never $120^\circ$.
Define a {\dfn gentle loop} as a gentle path such that the first and
last edges coincide. The smallest puzzles with gentle loops are of
size $6$; the one in figure \ref{fig:puzregions} is one of the only
two of that size.

(For this definition we did not really need to introduce puzzle regions, 
only region edges. We will need the regions themselves in 
section \ref{sec:puzresults}.)

\begin{figure}[htbp]
  \begin{center}
    \leavevmode
    \epsfig{file=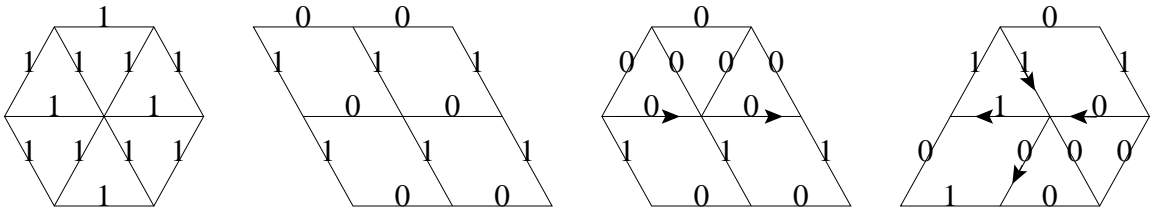,width=6in}
    \caption{The possible local structures of a puzzle near an
      interior vertex, up to rotation and puzzle duality. 
      The region edges incident on the vertex are oriented.
      This makes visible a ``traffic-planning'' mnemonic about gentle turns:
      any turn off a straightaway onto side roads is gentle, 
      whereas there are no gentle on-ramps onto straightaways.}
    \label{fig:puzlocal}
  \end{center}
\end{figure}

To better understand gentle paths, we need to know the possible local
structures of a puzzle around an interior vertex, which are
straightforward to enumerate. Clockwise around a lattice point in a
puzzle, we meet one of the following (see figure \ref{fig:puzlocal} 
for examples of each):
\begin{itemize}
\item six triangles of the same type
\item four rhombi at acute, obtuse, acute, obtuse vertices
\item three triangles of the same type, then two rhombi
\item some $0$-triangles, an acute rhombus vertex, some $1$-triangles,
  and an obtuse rhombus vertex.
\end{itemize}
Only the latter two have region edges.
This fourth type we call a {\dfn rake} vertex of the puzzle. 
(The reason for the terminology will become
clearer in lemma \ref{lem:rakes}.)

From this list, one sees easily that puzzle regions are necessarily
convex: traversing their boundaries clockwise, one never turns left.
In particular rhombus regions are necessarily parallelograms.

\subsection{Puzzles have either witnesses or gentle loops.}\label{ss:notaries}
Let $P$ be a puzzle constructed from a clockwise overlay $A\ovl B$ 
of two generic honeycombs, as in lemma \ref{lem:witness}.
A gentle path in $P$ is in particular a sequence of puzzle edges, each
successive pair sharing a vertex; there is a corresponding sequence
of edges in $A\ovl B$, each successive pair being sides of the same region.

\begin{Proposition}\label{prop:gpaths}
  Let $A\ovl B$ be a clockwise overlay of two generic honeycombs,
  $P$ the corresponding puzzle (as in theorem \ref{thm:klypuz}),
  $\gamma = (\gamma_1,\ldots,\gamma_m)$ a gentle path in $P$, 
  and $\tilde\gamma$ the corresponding sequence of edges in $A\ovl B$.
  Then the edge $\tilde\gamma_1$ is strictly longer than the edge
  $\tilde\gamma_m$.
\end{Proposition}

\begin{proof}
  It is enough to prove it for gentle paths of length two, and then
  string the $m-1$ inequalities together. In figure \ref{fig:gturns}
  we present all length two paths (up to rotation and dualization), 
  and the corresponding pairs of edges in an overlay. In each case the
  angles around the associated region force the strict inequality.

\begin{figure}[htbp]
  \begin{center}
    \leavevmode
    \epsfig{file=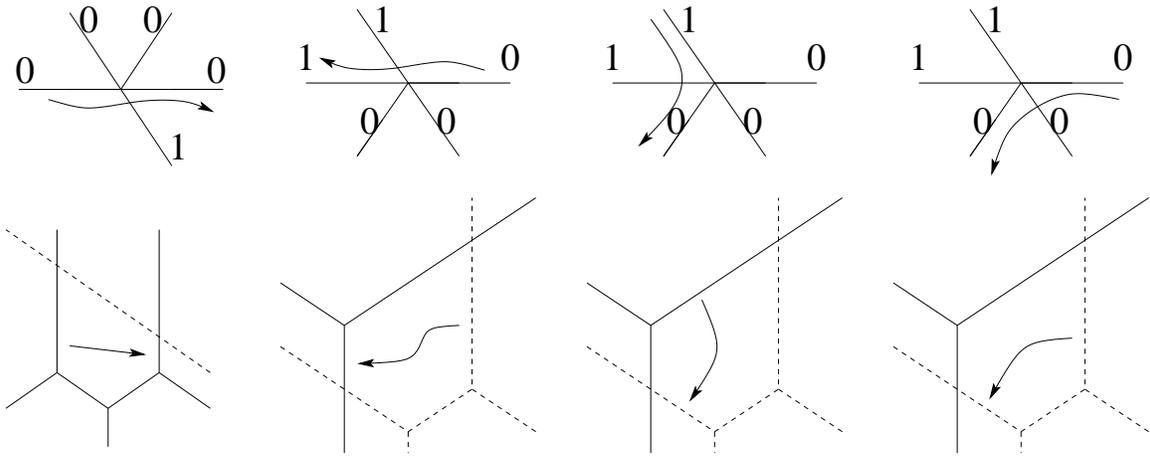,width=6in}
    \caption{Gentle turns (through the latter two diagrams in 
      figure \ref{fig:puzlocal}) and the associated pairs of edges in
      a generic witness, up to rotation and dualization.  Note that in
      each of the four gentle paths shown, the initial edge is
      necessarily longer than the terminal edge.  Dualization flips the
      pictures, and exchanges $0$ for $1$, giving the same geometric
      inequality.}
    \label{fig:gturns}
  \end{center}
\end{figure}
\end{proof}

\begin{Corollary}
  Let $P$ be a puzzle with a gentle loop.  Then $P$ does not arise
  from a clockwise overlay (``no witnesses''), and the inequality $P$
  gives on $\LRn$ is inessential.
\end{Corollary}

\begin{proof}
  If $P$ arises from a clockwise overlay $A\ovl B$, we can perturb
  $A$ and $B$ a bit to make them generic (as in the proof of
  theorem \ref{thm:klypuz}). Recall that a gentle loop is a gentle path 
  whose first and last edges agree. Then by the proposition the 
  corresponding edge in $A\ovl B$ is strictly longer than itself,
  contradiction.

  By the contrapositive of lemma \ref{lem:llco}, 
  $P$'s inequality is inessential.
\end{proof}

We now show that, conversely, gentle loops are the \emph{only}
obstructions to having witnesses.  In other words, if a puzzle
contains no gentle loops, then one can construct a witness $h$.  
The coming proposition \ref{prop:notaries} is inspired by the Wiener
path integral, in which a solution to a PDE at a point $x$ is
constructed as the sum of some functional over all possible paths from
$x$ to the boundary.  In our situation the role of the PDE is played
by the requirement that the edges around a region of $h$ close up to
form a polygon.  This construction will give a witness to the puzzle
provided that the number of gentle paths is finite, or equivalently if
there are no gentle loops.

\begin{Lemma}\label{lem:rakes}
  Let $P$ be a puzzle without gentle loops, and $v$ a rake vertex,
  as in figure \ref{fig:rake}, where four region edges meet. 
  Call the east edge on $v$ the {\dfn handle} and the west edges the
  {\dfn tines}. Then the number of gentle paths 
  starting at each of those four edges and terminating at
  the boundary is $a+b$ from the handle,
  $a$ from the two outer tines, and $b$ from the inner tine, for some
  $a,b\in \naturals$. 

\begin{figure}[htbp]
  \begin{center}
    \leavevmode
    \epsfig{file=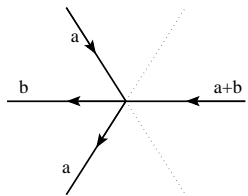,height=1in}
    \caption{Every vertex of a puzzle region, interior to a puzzle,
      looks like a rake (up to rotation and reflection). This one
      has the ``handle'' on the right, and the three ``tines'' on the left.}
    \label{fig:rake}
  \end{center}
\end{figure}

\end{Lemma}

\begin{proof}
%  These numbers are finite exactly because $P$ has no gentle loops.

  Of the two outer tines of the rake, one points toward the vertex, 
  one points away.
  Every gentle path starting at the inward-pointing tine goes into the
  outward-pointing tine (to turn into the middle tine would not be gentle), 
  and conversely every path from the outward-pointing can be extended;
  this is why they both have $a$ gentle paths to them for some 
  $a\in \naturals$.
  
  The gentle paths starting at the handle go either into the middle tine, 
  or the outward-pointing outer tine. So if $b\in \naturals$ gentle paths
  start from the middle tine, $a+b$ start from the handle.
\end{proof}

We will use an equivalent geometrical statement: these numbers
$a,b,a,a+b$ form the side lengths of a trapezoid in the triangular lattice.

\begin{Proposition}\label{prop:notaries}
  Let $P$ be a puzzle of size $n$ with no gentle loops. Then there
  exists a clockwise overlay $h = A\ovl B$ such that 
  the puzzle that theorem \ref{thm:klypuz} associates to $h$ is $P$.
  Therefore by lemma \ref{lem:witness}, the inequality defined by
  the puzzle determines a facet of $\LRn$.

  Moreover, in this witness $h$ every bounded nonzero edge of $A$ 
  crosses a bounded nonzero edge of $B$, and vice versa.
\end{Proposition}

We do this by direct construction. One example to follow along with
is the left honeycomb in figure \ref{fig:tco2puz}, whose regions are
indeed all trapezoids.

\begin{proof}
  First note that for this purpose, it is enough to specify a
  honeycomb {\em up to translation.}
  To do that, it is enough to specify the lengths and multiplicities
  of the bounded edges, and say how to connect them. Not every specification
  works -- the vertices have to satisfy the zero-tension condition
  of \cite{Hon1}, and %(dually) 
  the vector sum of the edges around a region must be zero.

  To specify a clockwise overlay, we must in addition two-color the edges
  ``$A$'' and ``$B$'', and make such edges only meet at crossing vertices,
  such that the colors alternate $A,B,A,B$ when read clockwise around
  each crossing vertex.

  We now build a clockwise overlay $h$ as a sort of 
  graph-theoretic dual\footnote{{\em} not in the sense of ``puzzle duality''}
  of $P$ -- one vertex for each puzzle region of $P$, one bounded edge
  for each region edge of $P$ (almost -- two region edges on the same
  boundary of a region determine the same edge of $P$), one semiinfinite 
  edge for each exterior edge of $P$.  If $e_h$ is a bounded edge of
  $h$ corresponding to a region edge $e_P$ of $P$, then $e_h$
  \begin{itemize}
  \item is perpendicular to $e_P$
  \item is labeled $A$ or $B$ depending on whether $e_P$ is adjacent
    to a $1$ or $0$ region
  \item has multiplicity equal to the length of $e_P$
  \item has length equal to the number of gentle paths starting
    at $e_P$ and ending on the boundary of $P$.
  \end{itemize}
  This last number is finite exactly because there are no gentle loops.

  The vertices of $h$ are zero-tension because 
  the vector sum of the edges around a region in $P$ is zero. 
  The regions in $h$ are all trapezoids, dual to the meeting of
  four regions in $P$ at rakes, and close up by lemma \ref{lem:rakes}.
  
  For the second statement, note that no edge of $h$ connects two distinct
  vertices in the same honeycomb, because the corresponding region edge
  in $P$ does not bound two regions of the same type.
\end{proof}

%A slight variant of this construction produces {\em all} the witnesses
%to $P$ satisfying the second condition of the proposition; put 
%positive real labels on the boundary vertices of the puzzle, and
%weight the gentle paths by their initial point.
%
The witnesses produced by this construction seem so minimal and
natural that we are tempted to christen them ``notaries''.
It would be interesting if there are correspondingly canonical 
witnesses in the Hermitian matrix context.

Together, this proposition \ref{prop:notaries} and the corollary to
proposition \ref{prop:gpaths} prove

\begin{Theorem}\label{thm:gloops}
  There is a $1:1$ correspondence between $n$-puzzles without gentle
  loops and regular facets of $\LRn$, given by the assignment
  $P\mapsto ``f_P \geq 0''$.
\end{Theorem}

At this point we have a complete, combinatorial characterization of 
the regular facets of $\LRn$ -- they correspond one-to-one to puzzles of
size $n$ with no gentle loops. However, to better tie in to Belkale's
result we need to characterize such puzzles in terms of rigidity.
One direction (Belkale's) is the following theorem, the other to come in the
next subsection. 

\begin{Theorem}\label{thm:belkpuz}
  Let $P$ be a puzzle with no gentle loops. Then $P$ is rigid.  
  In particular (by theorem \ref{thm:gloops}), 
  the set of rigid puzzles of size $n$ gives a 
  complete set of inequalities determining $\LRn$.
\end{Theorem}

We use in this proof one result that does not come until proposition
\ref{prop:puzstructure}: the number of rhombi in a puzzle is 
determinable from the boundary conditions.

\begin{proof}
  Let $h = A\ovl B$ be a witness to $P$ produced during the proof of
  proposition \ref{prop:notaries}. By slight perturbation of $A$ and
  $B$ to nondegenerate $A'$ and $B'$, small enough that no vertex of
  one crosses an edge of the other, we can create an $h'$ whose only
  degenerate edges correspond to the rhombi in $P$. 
  (One can modify the construction in proposition \ref{prop:notaries}
  to give such an $h'$ directly, but it is not especially enlightening.)
  
  Puzzles are easily seen to be determined by the set of their rhombi,
  and by proposition \ref{prop:puzstructure} (to come) the number of
  rhombi in a puzzle is determinable from the boundary conditions. So
  if $P$ is not rigid, so there exists another puzzle $Q$ with the
  same boundary, then this $Q$ has a rhombus that $P$ doesn't. Then by
  its definition as a sum of edge-lengths, $f_Q(h') > 0$. But by
  theorem \ref{thm:hrpuz}, $f_Q(h') = f_P(h') = 0$, contradiction.
\end{proof}

One way to think about this is that if a puzzle inequality can be
``overproved,'' there being two distinct puzzles $P,Q$ giving the same
inequality $f_P = f_Q \geq 0$, then the inequality is inessential.

(We will prove a stronger version of this result, in theorem 
\ref{thm:belkstrong}.)

\subsection{Breathing gentle loops.}\label{ss:breathing}

It remains to be shown that puzzles with gentle loops are not rigid.
The proof of this is very direct; given a sufficiently nice gentle loop 
$\gamma$ in a puzzle $P$, we will modify $P$ in a radius-$1$ neighborhood
of $\gamma$ to get a new puzzle $P'$ agreeing with $P$ outside that
neighborhood, in particular on the boundary. 
% (An example is in figure \ref{fig:breathex}.) 
The technical part comes in showing that minimal gentle loops are
``sufficiently nice.''

Define the {\dfn normal line} to a vertex $v$ along a 
gentle path $\gamma$ to be a pair of edges attached to $v$ such that
\begin{itemize}
\item they are $180^\circ$ apart
\item neither is in $\gamma$
\item neither cuts through the middle of a rhombus puzzle piece.
\end{itemize}
Checking the four cases in figure \ref{fig:gturns}, one sees that 
a normal line exists uniquely at each $v$. Note that the half of the
normal line connected to the left side of $\gamma$ is always labeled
$0$, and the right half always $1$.

We haven't needed to speak of the {\dfn distance between two puzzle vertices}
before; define it to be the graph-theoretic distance, where the graph
in question is made from the lattice triangle's vertices and edges
({\em not} just the edges appearing in the puzzle).

\begin{Lemma}\label{lem:loopbreathing}
  Let $P$ be a puzzle with a gentle loop $\gamma$, such that the only
  pairs of $\gamma$-vertices that are at distance $1$ in
  the puzzle are consecutive in $\gamma$. (In particular, the loop
  does not cross itself.) Then there exists a different puzzle 
  $P' \neq P$ that agrees with $P$ on any edge not touching $\gamma$.
\end{Lemma}

\begin{proof}
  Let $N\gamma$ denote the radius-$1$ neighborhood of $\gamma$,
  i.e. the set of pieces of $P$ with a vertex on $\gamma$.
  By the condition about nonconsecutive vertices, this neighborhood
  doesn't overlap itself, i.e. every edge connected to $\gamma$ (but
  not in $\gamma$) is connected to a {\em unique} vertex of $\gamma$.
  
  Cut $N\gamma$ up along its normal lines. It is easy to check that it
  falls into only four kinds of ``assemblages'' up to rotation, listed
  in figure \ref{fig:assemblages}.  
  \begin{figure}[htbp]
    \begin{center}
      \leavevmode
      \epsfig{file=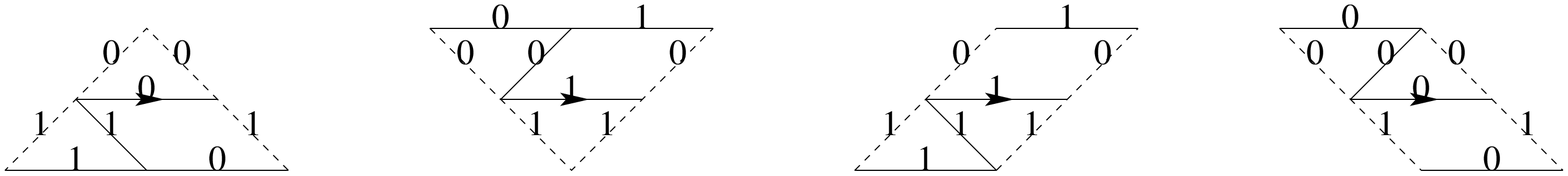,width=6in,height=1.2in}
      \caption{The four assemblages possible when a neighborhood of 
        a gentle path is cut along normal lines (here drawn dashed).
        These have been rotated to make the gentle path point East.}
      \label{fig:assemblages}
    \end{center}
  \end{figure}
  Notice that for each assemblage, there is a unique other of 
  the same shape, but rotated $180^\circ$.  This pairs up the two
  triangular assemblages, the parallelograms each being self-matched.
  The crucial observation to make is that two matching assemblages 
  have {\em the same labels on the boundary} (away from the normal lines).

  In particular, if we simultaneously replace each assemblage in
  $N\gamma$ by the other one with the same shape, the new collection
  fits together (because all the normal lines have been reversed),
  fits into the rest of the original puzzle (because the labels on
  the boundary are the same), and gives a new gentle loop running
  in the opposite direction.
\end{proof}

\begin{figure}[htbp]
  \begin{center}
    \leavevmode
    \epsfig{file=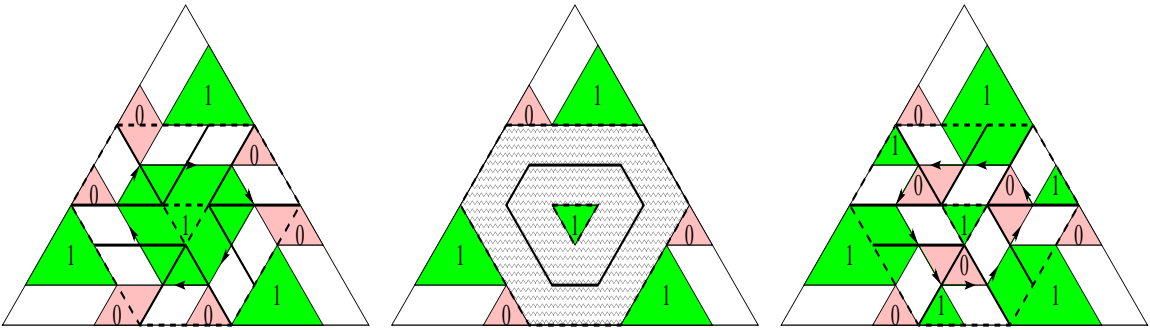,width=6in,height=1.7in}
    \caption{The left figure has a gentle loop running clockwise, 
      with its neighborhood broken along normal lines into
      assemblages, as indicated.
      Removing that neighborhood gives the middle figure, and
      filling it in with each assemblage replaced by its match gives the
      right figure, which has a gentle loop running counterclockwise.}
    \label{fig:breatheglobal}
  \end{center}
\end{figure}

We call this operation {\dfn breathing} the gentle loop, for
reasons explained at the end of section \ref{sec:puzresults}.
An example is given in figure \ref{fig:breatheglobal}.

(In fact the new puzzle constructed this way is unique, and breathing
the new gentle loop reproduces the original puzzle.)

\begin{Theorem}\label{thm:rigid}
  If a puzzle has gentle loops, it is not rigid.
\end{Theorem}

We will cut down the cases considered in this theorem using the
following lemma, easily checked from figure \ref{fig:gturns}:

\begin{Lemma}\label{lem:turns}
  If a 2-step gentle path $(\gamma_1,\gamma_2)$ in a puzzle $P$ turns 
  while passing through a vertex (as opposed to going straight), 
  there is another path $(\gamma_1',\gamma_2)$
  turning from the opposite direction.
  %, and in between them is a rhombus piece.
\end{Lemma}

\begin{proof}[Proof of theorem \ref{thm:rigid}]
  We will show that {\em minimal} gentle loops
  satisfy the condition of lemma \ref{lem:loopbreathing}.
  
  First we claim that minimal gentle loops do not self-intersect.
  Let $\gamma$ be a gentle loop that does self-intersect, and let $v$
  be a vertex occurring twice on $\gamma$ such that the two routes
  through $v$ are different. (If there is no such $v$, then $\gamma$
  is just a %gratuitously
  repeated traversal of a loop that does not self-intersect.) There
  are ${3\choose 2}=3$ local possibilities, corresponding to choosing
  two of the rightmost three gentle paths in figure
  \ref{fig:gturns}; in each one we can break and reconnect the
  gentle loop to make a shorter one, contradicting minimality.
  
  Second (and this is the rest of the proof) we claim that minimal
  gentle loops do not have nonconsecutive vertices at distance $1$. 
  If $\gamma$ is a counterexample, then there exists an edge 
  $E \notin \gamma$ connecting two points on $\gamma$. (This $E$ is
  just an edge in the lattice, not necessarily in the puzzle $P$ -- 
  it may bisect a rhombus of $P$ or whatever.)
  % ; call this a {\dfn bottleneck} of $\gamma$.
  Removing the endpoints of $E$ from $\gamma$ separates $\gamma$ into
  two arcs; call the shorter one the {\dfn minor arc} and the longer the
  {\dfn major arc}. (If they are the same length make the choice 
  arbitrarily.) Choose $E$ such that the minor arc is of minimal length.

  For the remainder we assume (using puzzle duality if necessary)
  that the gentle loop is clockwise.
%, the
%  counterclockwise case being similar. (Note that the two cases
%  are {\em not} exchanged by .)
  We now analyze the local picture near $E$, which
  for purposes of discussion we rotate to horizontal so that the minor
  arc starts at the west vertex of $E$, and ends at the east vertex.
  This analysis proceeds by a series of reductions, pictured in figure
  \ref{fig:reductions}.
  \begin{figure}[htbp]
    \begin{center}
      \leavevmode
      \epsfig{file=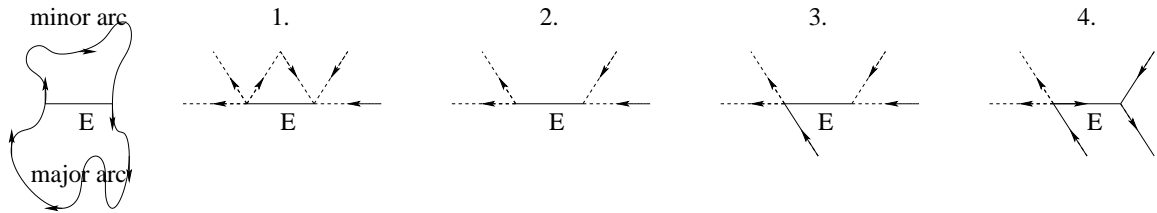,width=6in}
      \caption{The steps of theorem \ref{thm:rigid}. The dashed lines
        indicate possible edges at each step 1-4, the solid lines
        definite ones.}
      \label{fig:reductions}
    \end{center}
  \end{figure}

  1. The first edge of the minor arc goes either west, northwest, or northeast
  (to have room for a gentle turn from the major arc);
  likewise, the last edge goes either southeast, southwest, or west. 
  
  2. If the first edge of the minor arc went northeast, we could shift
  $E$ to $E'$ connecting the second vertex of the minor arc to the
  last vertex, contradicting the assumption that the minor arc was
  minimal length. So in fact the first edge goes west or northwest,
  the last southwest or west (by the symmetric argument).
  
  3. We now involve the major arc. Its last edge % of the major arc
  goes northwest or northeast.  If it goes northeast, then the first
  edge of the minor arc goes northwest (for gentleness). But then by lemma
  \ref{lem:turns} $E$ is oriented west; therefore we could shorten
  $\gamma$ to a loop that used $E$, contradicting $\gamma$'s assumed
  minimality. So the last edge of the major arc goes northwest.
  
  4. Therefore the first edge of the major arc goes southeast (since
  $\gamma$ doesn't intersect itself), so the last edge of the minor
  arc goes southwest. But then lemma \ref{lem:turns} says that $E$ is
  oriented east.
  
  At this point the west vertex of $E$ has 
  an oriented edge coming in from the southeast, one going out 
  to the east, and another going out either west or northwest
  (at least). This matches none of the vertices (or their puzzle duals) 
  in figure \ref{fig:puzlocal}.
  The contradiction is complete; there was no such $E$.
%  The case of $\gamma$ counterclockwise is exactly the same, exchanging
%  ``east'' and ``west'' everywhere in the argument.
\end{proof}

The following strengthening of theorem \ref{thm:belkpuz}
was also observed experimentally by W. Fulton.

\begin{Theorem}\label{thm:belkstrong}
  Let $P$ be a nonrigid puzzle. Then the face of $\LRn$ 
  determined by $f_P = 0$ lies on a chamber wall.
\end{Theorem}

\begin{proof}
  If not, there exists a regular boundary $b\in \LRn$ such that $f_P(b) = 0$.
  Let $h$ be a largest lift of $b$; by theorem 2 of \cite{Hon1} $h$
  is simply degenerate. Since $P$ is nonrigid, by theorem
  \ref{thm:belkpuz} it has a gentle loop; a minimal such loop
  $\gamma$ is breathable, by the proof of theorem \ref{thm:rigid}.
  Let $P'$ be the result of breathing $P$ along $\gamma$.
  
  We claim that some rhombus $\rho'$ of $P'$ overlaps some rhombus
  $\rho$ of $P$ in a triangle. To see this, divide $\gamma$ up along
  its normal lines into the assemblages of lemma
  \ref{lem:loopbreathing}, figure \ref{fig:assemblages}.  At least
  one such assemblage must be a triangle, for otherwise the loop can
  not close up. (In fact there must be at least six triangles.) When
  we breathe the loop, the rhombus $\rho$ in the original assemblage in $P$
  overlaps the rhombus $\rho'$ in the new assemblage in $P'$ in a triangle.
  
  Dually, the edges of the honeycomb tinkertoy $\tau_n$ corresponding
  to those two rhombi meet at a vertex.

  Since $f_P(h) = f_P'(h) = 0$, and they are defined as the sum of 
  certain edge-lengths of $h$, all those edges of $h$ must be length zero.
  Therefore, the two {\em adjacent} edges in $h$ corresponding to
  $\rho$ and $\rho'$ are length zero.
  But then $h$ is not simply degenerate, contradiction.
\end{proof}

At this point we have a second characterization of the regular facets
of $\LRn$; they correspond to rigid puzzles. For our final characterization
we need to involve the other life of puzzles, which is in computing
Schubert calculus of Grassmannians.

%Non-honeycomb-based proofs of this result are also known.

\section{Puzzle inflation, rhombi, and Schubert calculus}\label{sec:puzresults}

We start with an ``inflation'' operation on puzzles, taking a puzzle
$P$ and a natural number $N$ to a new puzzle $N\cdot P$:

\begin{Lemma}\label{lem:Npuz}
  Let $P$ be a puzzle of size $n$. For $N\in \naturals$, define
  $N\cdot P$ to be the puzzle whose puzzle regions
  are in correspondence with $P$'s, and glued together the same way,
  but every edge labeled $1$ has been stretched by the factor $N$.
  Then $N\cdot P$ is a well-defined puzzle. 
  In addition, $N\cdot P$ is rigid if and only if $P$ is rigid.
\end{Lemma}

\begin{figure}[htbp]
  \begin{center}
    \leavevmode
    \epsfig{file=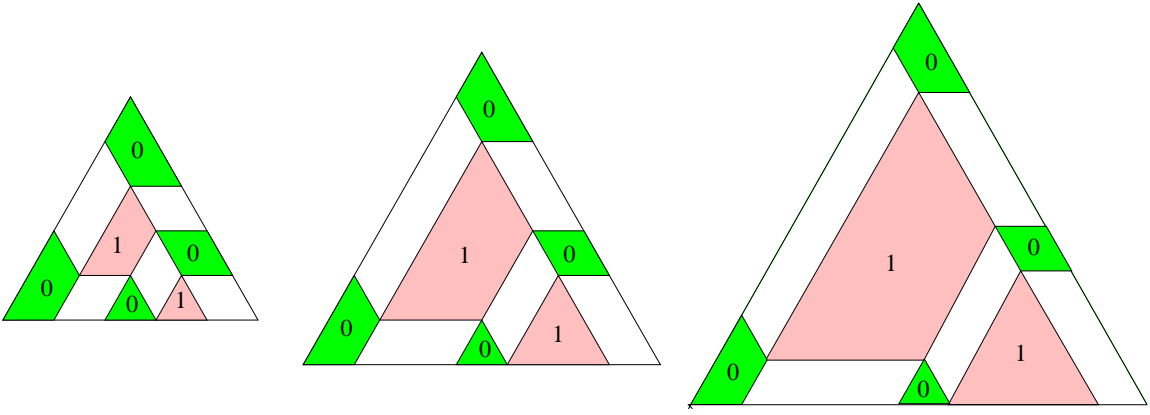,width=6in}
    \caption{A puzzle $P$ broken into regions, $2\cdot P$, and $3\cdot P$.}
    \label{fig:Npuz}
  \end{center}
\end{figure}

An example is in figure \ref{fig:Npuz}.

\begin{proof}
  To see that $N\cdot P$ is a well-defined puzzle,
  we need to check that each new region is well-defined -- that
  traversing its boundary we return to where we started.
  For $0$-regions there is nothing to show, for $1$-regions the whole
  region is inflated by the factor $N$, and for rhombus regions
  two opposite sides of the parallelogram are stretched.

  There is an evident correspondence between gentle loops in $P$
  and gentle loops in $N\cdot P$, so by theorems \ref{thm:belkpuz} 
  and \ref{thm:rigid} they are either both rigid or neither is.
\end{proof}

Our principal use of this will be the {\dfn deflation} $0\cdot P$
of a puzzle.

\begin{Proposition}\label{prop:puzstructure}
  Let $P$ be an $n$-puzzle, such that one side has $r$ $1$-edges.
  Then the other two sides also have $r$ $1$-edges, and the puzzle
  consists of 
  \begin{itemize}
  \item $r^2$ $1$-triangles, of which ${r+1 \choose 2}$ are right-side-up and
    ${r \choose 2}$ are upside-down
  \item $(n-r)^2$ $0$-triangles, of which ${(n-r)+1 \choose 2}$ are
    right-side-up and ${n-r \choose 2}$ are upside-down
  \item $r(n-r)$ rhombi.
  \end{itemize}
\end{Proposition}

Since we already know every facet comes from a puzzle 
(theorem \ref{thm:klypuz}), this implies Horn's results on the structure
of the inequalities determining facets, explained in 
subsection \ref{ss:horn}.

\begin{proof}
  The deflation of $P$ is an $(n-r)$-puzzle consisting only of
  $0$-triangles, letting us count the number of $0$-triangles in the
  original puzzle (namely, $(n-r)^2$) and also the number of $0$-edges
  on the three sides (namely, $n-r$). By deflating the dual puzzle,
  one can count there to be $r^2$ $1$-triangles.
  The remaining area is $n^2 - r^2 - (n-r)^2 = 2r(n-r)$ in units of
  triangle, and all that is available are rhombi which each use up 2.
\end{proof}

We now prove that puzzles compute Schubert calculus 
on Grassmannians (our reference for the latter is \cite{F2}),
as stated in the introduction. 
Recall that we index Schubert classes in $H^*(\Grrn)$ by $n$-tuples
consisting of $r$ ones and $n-r$ zeroes.\footnote{
It is more usual to encode these classes by partitions fitting inside
an $(n-r) \times r$ rectangle. The correspondence is as follows.
Given one of our $n$-tuples, read the $0$s
as ``left'' and the $1$s as ``down''; this gives a path from the upper right
corner of such a rectangle to the lower left.
Above this is a partition, and conversely, given a partition in this
rectangle we can read off an $n$-tuple of ``left''s and ``down''s.
}

What we will actually prove is that the puzzle rule is equivalent to
the honeycomb rule from \cite{Hon1}. (A direct proof will appear in
\cite{Puz1}, in turn giving an independent proof of the honeycomb rule.)
First we need a lemma on honeycombs:

\begin{Lemma}\label{lem:blastradius}
  Let $h$ be a honeycomb with boundary coordinates 
  $(\lambda,\mu,\nu) \in (\reals^n)^3$ on the 
  Northwest, Northeast, and South sides, 
  each a weakly decreasing list of real numbers. Then
  \begin{enumerate}
  \item The first coordinate of any vertex of $h$ is in the interval
    $[\lambda_n,\lambda_1]$. (Likewise second coordinate in $[\mu_n,\mu_1]$, 
    third coordinate in $[\nu_n,\nu_1]$.)
  \item The third coordinate of any vertex is in the interval
    $[-\lambda_1-\mu_1, -\lambda_n-\mu_n]$.
  \item If $\lambda_i,\mu_i,-\nu_i \in [0,M]$ for all $i=1\ldots n$,
    then all of $h$'s vertices are in the triangle with vertices 
    $(0,0,0), (M,0,-M), (0,M,-M)$.
  \end{enumerate}
\end{Lemma}

\begin{proof}
  1. Follow a path in the honeycomb, going Northwest whenever possible,
  Southwest when not, eventually coming out on an edge with constant
  coordinate $\lambda_i \leq \lambda_1$. Each Southwest sojourn increases
  the first coordinate, and each Northwest leaves it unchanged, so the
  original first coordinate must have been at most $\lambda_1$.
  Replacing ``Southwest'' with ``North'' gives the opposite inequality.
  The other two coordinates come from rotating this proof $\pm 120^\circ$.
  
  2. Since the sum of the three coordinates is zero by definition, the
  third one can be bounded in terms of the first two.

  3. This is just a special case of (1).
\end{proof}

\begin{Theorem*}[theorem \ref{thm:puzcount} from the introduction]
  Let $\pi,\rho,\sigma$ be three strings of $r$ ones and $n-r$ zeroes.
  Then the number of puzzles with $\pi,\rho,\sigma$ clockwise 
  around the boundary is the Schubert intersection number 
  $\int_{\Grrn} S_\pi S_\rho S_\sigma$.

  Equivalently, write $S_\pi S_\rho = \sum_\tau c_{\pi\rho}^\tau S_\tau$.
  The the number of puzzles with $\pi,\rho$ on the NW and NE
  boundaries, and $\tau$ on the South boundary, all written left to right,
  is the structure constant $c_{\pi\rho}^\tau$.
\end{Theorem*}

\begin{proof}
  We prove the second statement: the first follows from the second,
  since $\int_{\Grrn} S_\pi S_\tau$ is $1$ if $\pi$ is the reversal of $\tau$,
  $0$ otherwise.

  The structure constants for multiplication of Schubert classes are
  well known to also be the structure constants for tensor products of
  polynomial representations of $\GLn$ (the first to observe this seems
  to be Ehresmann; see \cite{F2} or \cite{Puz2}).
  The precise statement is as follows. Let $\lambda_i$ be the number of
  $0$s after the $i$th $1$ in $\pi$, 
  so $\lambda_1 \geq \lambda_2 \geq \ldots \geq \lambda_r$, and 
  $\lambda$ is a partition of the number of inversions of $\pi$.
  Likewise construct $\mu$ from $\rho$, and $\nu$ from $\tau$.
  Then
  $$ c_{\pi\rho}^\tau = \dim \Hom_{\GLn}(V_\nu, V_\lambda \tensor V_\mu). $$
  This latter can be calculated using honeycombs, as proved in
  \cite{Hon1}; it is the number of honeycombs with boundary
  coordinates $(\lambda_1 \geq \ldots \geq\lambda_r)$ on the Northwest
  side, $(\mu_1 \geq \ldots \geq \mu_r)$ on the Northeast, and
  $(-\nu_r \geq \ldots \geq -\nu_1)$ on the South.

  We now construct a map from our puzzles to these honeycombs.
  To create a honeycomb from a puzzle is a three-step process
  (follow along with the example in figure \ref{fig:fulldeflate}):
  \begin{enumerate}
  \item Place the puzzle in the plane $\reals^3_{\sum=0}$ 
    such that the bottom right corner is at the origin, 
    and turn it $30^\circ$ counterclockwise.
  \item At each boundary edge labeled $1$, attach a rhombus 
    ({\em outside} the puzzle), then another (parallel to the first),
    and repeat forever. Fill in the rest of the plane with $0$-triangles.
  \item Deflate the extended puzzle, keeping the right corner at the origin. 
    The honeycomb's vertices then come from the deflated $1$-regions,
    and the honeycomb's edges come from the deflated rhombus regions,
    with the multiplicity on the edge coming from the thickness of the
    original rhombus region.
  \end{enumerate}

\begin{figure}[htbp]
  \begin{center}
    \leavevmode
    \epsfig{file=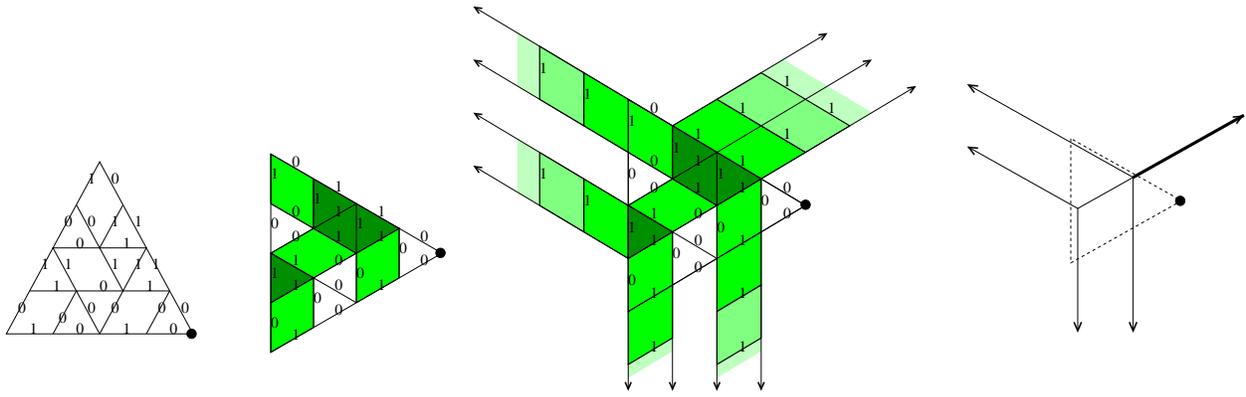,width=6.5in}
    \caption{A puzzle, and the three stages in creating a honeycomb from it.
      Regions with $1$s on their boundary are shaded, and then deflated.
      The origin $(0,0,0)$ is indicated in each figure by a heavy dot.}
    \label{fig:fulldeflate}
  \end{center}
\end{figure}

  The resulting diagram obviously has finitely many vertices, all edges in
  triangular-coordinate directions, and semiinfinite edges only going NW/NE/S
  (coming from the $1$-edges on the boundary of the original puzzle).
  The remaining condition for it to be the diagram of a honeycomb is
  that each vertex have zero total tension; this is equivalent to the fact
  that the original $1$-region was a closed polygon.

  To see that this is a bijection, we construct the inverse map.
  Start with a honeycomb computing $c_{\pi\rho}^\tau$. 
  By part 2 of lemma \ref{lem:blastradius}, it fits inside the triangle
  with vertices $(0,0,0), (n-r,0,r-n), (0,n-r,r-n)$.
  Inflate each edge of the honeycomb intersected with the triangle
  to a rhombus region, the thickness given by the multiplicity of the edge,
  and each vertex to a polygon of $1$-triangles, the lengths of the
  edges of the polygon given by the multiplicities of the edges 
  at the vertex. The result is a puzzle with boundary $\lambda,\mu,\nu$.
\end{proof}

In \cite{Puz1} will appear an alternate proof of this theorem not
using honeycombs, which shows also that puzzles compute 
{\em $T$-equivariant} cohomology of Grassmannians, when one includes
an additional ``equivariant puzzle piece''.

We conclude this section with some observations.

1. This intersection number problem has 12 manifest symmetries; 
$3!$ from permuting $\pi,\rho,\sigma$, and $2$ from the
duality diffeomorphism $\Gr_r(\complexes^n) \iso \Gr_{n-r}(\complexes^n)$.
As with honeycombs and Berenstein-Zelevinsky patterns, only half
of these are manifest in the rule; rotating the puzzles gives
the even permutations of $\pi,\rho,\sigma$, and puzzle duality
gives the composition of Grassmann duality with the odd permutation
$\pi \leftrightarrow \sigma$.
None of these are directly visible in the Littlewood-Richardson rule 
(for a deeper discussion of this, see \cite{Puz2}).

2. This theorem makes possible a (rather forced) duality on honeycombs,
as already observed in \cite{GP, Hon1}: pick a triangle containing the
honeycomb, inflate to a puzzle, apply puzzle duality, and deflate back
to a new honeycomb, the total effect being to exchange vertices for
regions and vice versa. Unfortunately this depends on the choice
of triangle, and only works for integral honeycombs. This is quite
different from the much more natural duality on honeycombs that comes
from flipping them over, $(x,y,z) \mapsto (-y,-x,-z)$.

3. In \cite{Hon1} we defined a way of locally modifying a honeycomb in
the vicinity of a loop through nondegenerate vertices, which was
also called breathing. It is easy to check that any breathing operation
on honeycombs is, under the deflation correspondence above, the deflation
of a gentle-loop breathing on a puzzle. (The reverse is not
true: gentle-loop breathing of puzzles is a strict generalization.)

4. Note that this connection of puzzles and honeycombs is completely
different from the one in theorem \ref{thm:klypuz}, and serves as a
combinatorial explanation of the recursive nature of Horn's list of
inequalities. To recapitulate the chain of reasoning involved: 
first one studies extremal $n$-honeycombs (as we did in section
\ref{sec:tco}), and from an extremal honeycomb, which is necessarily
an overlay of a $r$-honeycomb and an $(n-r)$-honeycomb, one constructs a
puzzle encoding the pattern of overlay. That puzzle then deflates to a
$r$-honeycomb, necessarily integral. Therefore inequalities on 
$n$-honeycombs can be ``blamed'' on integral $r$-honeycombs.
This is not recursive until one knows that honeycombs exist with given
integral boundary conditions if and only if integral honeycombs exist
with the same boundary;
this was (the honeycomb version of) the saturation conjecture, proved
in \cite{Hon1}.

\section{Replacing puzzles by Schubert calculus}\label{sec:schub}

So far we have used puzzles to give inequalities
on the boundaries of honeycombs. 
In this section we replace puzzles by Schubert calculus,
and honeycombs by zero-sum Hermitian triples,
to formulate puzzle- and honeycomb-free versions of most of our results.
The only casualty is the characterization of rigid puzzles as
those without gentle loops; but this too has an application,
in proving a conjecture of W. Fulton.

First we recall the connection of honeycombs to Hermitian matrices:
Let $(\lambda,\mu,\nu) \in (\reals^n)^3$ be three weakly decreasing
lists of real numbers. Then there exists a honeycomb with boundary
coordinates $(\lambda,\mu,\nu)$ if and only if there exists a
triple of Hermitian matrices $(H_\lambda,H_\mu,H_\nu)$ with spectra
$\lambda,\mu,\nu$ and adding to the zero matrix.
This follows from \cite{Hon1} and \cite{Kly}; we make a more precise
statement in the appendix, replacing Klyachko's argument with more direct
use of the relation between geometric invariant theory quotients and
symplectic quotients.

\begin{Corollary}
  Avoiding direct mention of honeycombs and puzzles, we have
  \begin{enumerate}
  \item \cite{T,HR,Kly}
    If $S_\pi,S_\rho,S_\sigma$ are three Schubert classes on $\Grrn$
    such that  $\int_{\Grrn} S_\pi S_\rho S_\sigma > 0$, then for
    any triple $H_\lambda,H_\mu,H_\nu$ of $n\times n$ Hermitian matrices
    with zero sum and spectra $\lambda,\mu,\nu$ (written in 
    decreasing order) respectively, we have the inequality
    $$ \sum_{i=1}^{n} \pi_i \lambda_i + \rho_i \mu_i + \sigma_i \nu_i \leq 0.$$
  \item \cite{Kly} This list of inequalities on the spectra
    is sufficient for the existence of such a triple.
  \item \cite{Be}
    If $S_\pi,S_\rho,S_\sigma$ are three Schubert classes on $\Grrn$
    such that  $\int_{\Grrn} S_\pi S_\rho S_\sigma > 1$, then the
    corresponding H-R/T/K inequality is inessential...
  \item 
   ...and equality can only occur when $(\lambda,\mu,\nu)$ are not all regular.
  \item 
    If $S_\pi,S_\rho,S_\sigma$ are three Schubert classes on $\Grrn$
    such that  $\int_{\Grrn} S_\pi S_\rho S_\sigma = 1$,  then the
    corresponding H-R/T/K inequality is essential.
  \end{enumerate}
\end{Corollary}

\begin{proof}
  These are theorem \ref{thm:puzcount} combined with
  \begin{enumerate}
  \item theorem \ref{thm:hrpuz} (the inequality is reversed because we are
    summing over $1$s here instead of $0$s)
  \item theorem \ref{thm:klypuz} 
  \item theorem \ref{thm:belkpuz} 
  \item theorem \ref{thm:belkstrong}
  \item theorems \ref{thm:gloops} and \ref{thm:rigid}.
  \end{enumerate}
\end{proof}

We have another puzzle-free application of theorem \ref{thm:puzcount}:

\subsection{Fulton's conjecture.}\label{ss:fulton}
In a private communication, W. Fulton proposed the following 

\begin{Conjecture}
Let $\lambda,\mu,\nu \in (\integers^n)^3$ 
be a triple of dominant weights for $\GLn$, 
and $V_\lambda,V_\mu,V_\nu$ the corresponding
irreducible representations.
% \hfill \break
If $V_\nu$ occurs exactly once as a constituent of $V_\lambda\tensor V_\mu$,
then $\forall N\in\naturals$, $V_{N\nu}$ occurs exactly once as 
a constituent of $V_{N\lambda}\tensor V_{N\mu}$.
\end{Conjecture}

It is interesting to compare this to the saturation conjecture (proven
in \cite{Hon1}). Saturation says that if a polytope of honeycombs with
fixed integral boundary is nonempty, the polytope contains at least
one lattice honeycomb.  The present conjecture is sort of a next step:
its contrapositive says that if a polytope of honeycombs with 
fixed integral boundary is not only nonempty but positive-dimensional,
the polytope contains at least {\em two} lattice honeycombs.

We need one additional construction in order to prove this conjecture:
the {\dfn dual inflation} of a puzzle by a factor $M$, defined as 
dualizing the puzzle, $M$-inflating, then dualizing again.
This amounts to thinking of the inflation of a puzzle in terms of the
$0$-edges instead of the $1$-edges.

\begin{proof}[Proof of Fulton's conjecture.]
  Let $\det = (1,1,\ldots,1)$ denote the high weight of 
  the determinant representation, and $c_{\lambda\mu}^\nu$ denote the
  number of times $V_\nu$ appears in $V_\lambda \tensor V_\mu$. 
  Then using the equality
$$ c_{\lambda\mu}^\nu = c_{\lambda+L\det,\,\mu+M\det}^{\nu+(L+M)\det}, $$
  we can reduce to the case that $\lambda$ and $\mu$ are nonnegative
  (so, high weights of polynomial representations).
  Therefore $\nu$ is also nonnegative, 
  for otherwise $c_{N\lambda,N\mu}^{N\nu}$ would be zero for all 
  $N\in \naturals$.

  From there, we can use theorem \ref{thm:puzcount} to
  convert to a Schubert problem, i.e. counting puzzles rather than honeycombs. 
  One then has to check that rescaling a honeycomb by the factor $N$
  corresponds to dual-inflation on puzzles.
  
  Since the original honeycomb $h$ is rigid, so too is the 
  corresponding puzzle $P$, therefore by lemma \ref{lem:Npuz} so too
  is the dual inflation of $P$ by the factor $N$, 
  and therefore so is $N\cdot h$.
\end{proof}

\subsection{$(N,M)$-inflation and non-polynomiality.}
Define the $(N,M)$-inflation of a puzzle by $N$-inflating it, and then
$M$-dual-inflating it (these operations commute).

Note that this descends to a well-defined notion of the $(N,M)$-inflation of a 
{\em boundary} condition on a puzzle, and as such one can study the functions
$$ f(N,M) := \# \hbox{puzzles with boundary 
        $(N,M)\cdot \pi, \, (N,M)\cdot \rho, \, (N,M)\cdot \sigma$} $$ 
for fixed initial boundary conditions $\pi,\rho,\sigma$.
Because of the connections of puzzles to honeycombs and thereby to
sections of a line bundle over $(\GLn/B)^3 // \GLn$ (see the appendix),
one can show that $f(N,M)$ is a polynomial function of one 
argument when the other is held fixed.\footnote{%
Geometrically, this is essentially due to the fact that the GIT quotients
$(\GLn/B)^3 // \GLn$ are usually {\em manifolds} and never {\em orbifolds},
a fact special to the group $\GLn$. A different proof is given in
\cite{DW2}.}

Taken together, though, the growth is usually exponential;
the reader may enjoy showing that for $\pi = \rho =\sigma = 010101$
(as in figure \ref{fig:puzregions}), $f(N,M) = {N+M \choose N}$.
%It seems a challenge to produce such a function in a way that generalizes
%the interpretation above of $f(N,1)$.

\section{Summing more than three matrices}

\newcommand\Bd[1]{\BDRY_{#1}(n)}
The cone $\LRn$, whose facets we have now completely determined, has
a generalization for any $m\in\naturals$: the set of $m$-tuples of
spectra 
$$ \Bd{m} := 
 \bigg\{ (\lambda_i)_{i=1\ldots m} \quad:\quad 
        \exists \{H_i\}, \sum_i H_i = 0 \bigg\}$$ 
such that there exist $n\times n$ Hermitian matrices with those spectra
adding to the zero matrix. (Again, this is equivalent to the 
corresponding $m$-fold tensor product problem.) Then $\Bd{3}$ is just the 
cone $\LRn$ we've already determined, and $\Bd{1},\Bd{2}$ are uninteresting.

To study this cone for $m>3$ by the techniques in \cite{Hon1} and this paper,
we need to determine the corresponding honeycomb extension problem. We do this
by factoring the problem: sum the first $m-2$ matrices and call the
eigenvalues of that $\mu$, then see if $\mu$ goes with the last two spectra.
$$ \Bd{m} = \bigg\{  (\lambda_i)_{i=1\ldots m} \quad :\quad \exists \mu,
(\lambda_1,\ldots,\lambda_{m-2},\mu^*) \in \Bd{m-1}, 
(\mu,\lambda_{m-1},\lambda_m) \in \Bd{3} \bigg\}$$
(Here $\mu^*$ denotes $-\mu$, reversed so as to again be in decreasing order.)
Repeat this factorization\footnote{%
In an alternate view of the Hermitian sum problem that we haven't 
discussed, about flat $U(n)$-connections on an $m$-punctured sphere with
small holonomies around the punctures,
this corresponds to taking a pants decomposition of the punctured surface.}
until everything is in terms of $\Bd{3} = \LRn$.
Then we can think of $\Bd{m}$ in terms of an $(m-2)$-tuple of
honeycombs such that one boundary of each honeycomb anti-agrees with
one boundary of the next.  Define an {\dfn $m$-ary honeycomb} as
exactly such an $(m-2)$-tuple.

\begin{figure}[Hhtb]
  \begin{center}
    \leavevmode
    \input{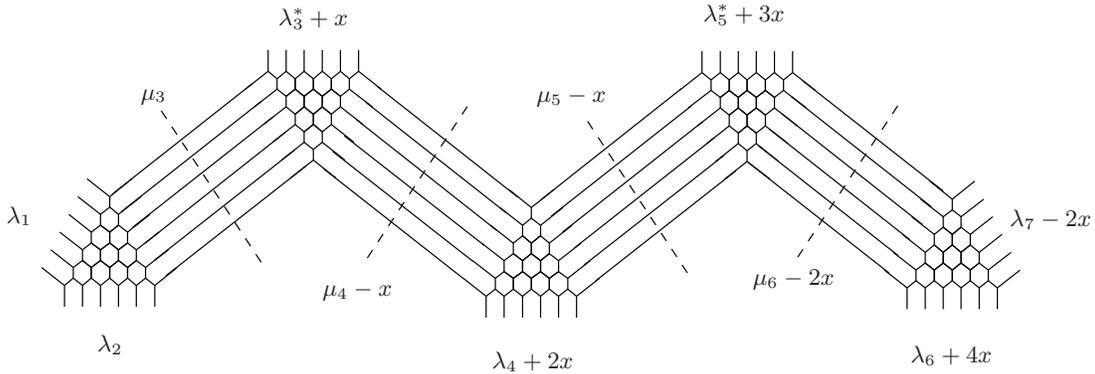}
    \caption{A honeycomb whose boundary lies in $\Bd{7}$.}
    \label{fig:multhoney}
  \end{center}
\end{figure}

Graphically, the easiest way to think about these is to draw the
honeycombs in the same plane, half of them upside down, 
and very far from one another, as in figure \ref{fig:multhoney}.
To get them far from one another we can add a large-enough constant $x$
to the coordinates.\footnote{%
  Mathematically, it is nicer to deal with the $(m-2)$-tuple, because
  it doesn't require one to choose this large-enough $x$. If one
  insists on actually working with these single composite diagrams,
  one must use part 3 of lemma \ref{lem:blastradius} bounding the
  size of a honeycomb in terms of two of its boundaries.}  
Note that we do {\em not} have to go beyond two dimensions, as is many
people's first guess about $m>3$ (or indeed $n>3$).

\begin{figure}[htbp]
  \begin{center}
    \leavevmode
    \epsfig{file=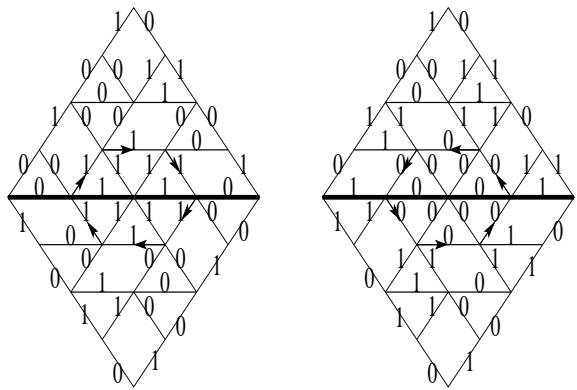,height=2in,width=3in}
    \caption{Given four generic lines in $\CP^3$, 
      exactly two lines touch all four. The gentle loops are drawn
      around the central hexagons.}
    \label{fig:4linestouch2}
  \end{center}
\end{figure}

All the same techniques developed for $m=3$ go through without change.
Define an {\dfn $m$-ary puzzle} as an $(m-2)$-tuple of puzzles (every
other one upside down) that can be fitted together into a line, and
call the individual puzzles the {\dfn constituents} of the $m$-ary puzzle.  
(Careful: these are not merely arrangements of puzzle pieces into a
trapezoid/parallelogram; they satisfy the extra condition that no
rhombus is allowed to cross from one constituent into the next.) 
On the Schubert calculus side, these count intersections of $m$ cycles in
a Grassmannian.  In figure \ref{fig:4linestouch2} we give the famous
count (two) of the number of lines touching four others in
$\CP^3$.

A {\dfn gentle path in an $m$-ary puzzle} has essentially the same definition,
with the only tricky point that it cannot include one of the edges
joining one constituent puzzle to the next. With these definitions we
have the analogous results:

\begin{Theorem}
  Each $m$-ary puzzle gives a nonnegative functional on $\Bd{m}$.  The
  regular facets of $\Bd{m}$ come precisely from the $m$-ary puzzles
  with no gentle loops. These $m$-ary puzzles are exactly the rigid ones,
  corresponding to $m$ Grassmannian Schubert cycles intersecting in
  a unique point.
\end{Theorem}

\begin{proof}
  All the proofs go through without modification, except for one:
  we need to check that when we breathe a gentle loop in an $m$-ary puzzle
  using the loop-breathing lemma \ref{lem:loopbreathing},
  we don't introduce any rhombi that cross from one constituent 
  puzzle to the next, for that would remove a boundary edge from a
  puzzle. But the edges separating constituents are obviously 
  on normal lines to the gentle path, and the loop-breathing construction
  does not remove these edges.
\end{proof}

\subsection{A representative example.}
In figure \ref{fig:4linestouch2} we exhibited two $4$-ary puzzles with
the same boundary.  By our theorems, we know that these have gentle
loops, and determine the same true inequality 
$$\lambda_1+\lambda_3 + \mu_1+\mu_3 + \nu_1+\nu_3 + \pi_1+\pi_3 \geq 0 $$
on spectra of four Hermitian matrices with zero sum,
but that this inequality is inessential. 

\section{Appendix: the equivalence of the definitions of $\LRn$}

All the arguments in this paper study $\LRn$ purely in terms of its
interpretation as the possible boundary conditions of honeycombs.
In \cite{Hon1}, these are related to invariants in tensor products
of $\GLn$-representations, which are turn related to Hermitian matrices
in \cite{Kly}.

We include here a stronger result
(which could already have been given in \cite{Hon1}),
replacing Klyachko's argument by the Kirwan/Kempf-Ness theorem,
allowing for a more precise result.
While this involves some somewhat formidable machinery, 
its application really is a routine matter, 
and so we label the following a corollary.

\begin{Corollary}[to theorem 4 of \cite{Hon1}]
Let $(\lambda,\mu,\nu) \in (\reals^n)^3$ be a triple of weakly
decreasing $n$-tuples of reals.
The volume (resp. real dimension) of the polytope of honeycombs 
in $\HONEY(\tau_n)$ with boundary coordinates $(\lambda,\mu,\nu)$
is equal to the symplectic volume (resp. complex dimension) 
of the space of zero-sum Hermitian triples with these spectra 
modulo the diagonal action of $U(n)$.
%If $(\lambda,\mu,\nu)$ are integral, this space can be identified
%with a moduli space of triples of flags in $\complexes^n$, and the
%volume of the polytope is the leading coefficient of the 
%Hilbert polynomial of this moduli space in its natural projective embedding.
%In particular this shows
%that the cone $\BDRY(\tau_n)$ is the same as the cone $LR_n$.
In particular, there exists such a honeycomb if and only if 
there exists such a zero-sum Hermitian triple.
\end{Corollary}

\begin{proof}
The machinery used here is geometric invariant theory,
particularly the ``geometric invariant theory quotients
are symplectic quotients'' theorem \cite{MFK}.\footnote{%
Klyachko's proof of the relation between these two problems
follows the same essential lines as this more general theorem.}
To begin with, take $(\lambda,\mu,\nu)$ integral, and consider the graded ring
$$       R := \bigoplus_k V_{k\lambda} \tensor V_{k\mu} \tensor V_{k\nu}. $$
By the Borel-Weil theorem, $\Proj R$ is a product of 
three (partial) flag manifolds as an algebraic variety, and from its
induced projective embedding inherits a symplectic structure.
By Kostant's extension of Borel-Weil, it is symplectomorphic to
the product of the $U(n)$ coadjoint orbits through
$\lambda,$ $\mu,$ and $\nu$.
We use the trace form on $\lie{u}(n)$ to identify these with
the corresponding isospectral sets of Hermitian matrices.

Now consider the invariant subring $R^{GL_n}$. By definition, 
$\Proj R^{GL_n}$ is the geometric invariant theory quotient of this
product of three flag manifolds by the diagonal action of $GL_n$.
By theorem 4 of \cite{Hon1}, the Hilbert function of this variety is
the Erhart function of the polytope of honeycombs with boundary
conditions $(\lambda,\mu,\nu)$. In particular its leading coefficient
(resp. its degree),
which as for any projective variety is the variety's symplectic volume
(resp. its complex dimension),
is also the volume (resp. real dimension) of the polytope.

By the GIT/symplectic equivalence, this GIT quotient by $GL_n$ can 
alternately be constructed as a symplectic quotient by its maximal 
compact $U(n)$. This construction takes the zero level set of 
the $U(n)$ moment map -- here the moment map is the sum of the three
matrices -- and quotients it by $U(n)$. Combining these results, 
we find that the symplectic volume (resp. complex dimension) of the 
moduli space of zero-sum Hermitian triples is the volume (resp. 
real dimension) of the polytope of honeycombs.

The same holds for rational triples $(\lambda,\mu,\nu)$ because both
sides behave the same way under rescaling, and then for arbitrary
triples because both sides are continuous.
\end{proof}

In particular (as Klyachko proves): an invariant tensor implies the
existence of a zero-sum Hermitian triple, and a zero-sum Hermitian
triple implies the existence of an invariant tensor in 
$V_{k\lambda} \tensor V_{k\mu} \tensor V_{k\nu}$ for some $k>0$.
To go from there to $k=1$ is harder, 
requiring the saturation conjecture proved in \cite{Hon1}.

%\goodpagebreak
\bibliographystyle{alpha}

\end{document}